\documentclass[final,pdftex]{siamltex}
\pdfoutput=1

\usepackage{amsmath,amssymb,amscd}
\usepackage{amsfonts}
\usepackage{graphicx}
\usepackage{grffile}
\usepackage[hang,centerlast]{subfigure}

\usepackage{tikz}

\newcommand{\ie}{i.e., }
\newcommand{\eg}{e.g., }
\newcommand{\lift}{\mathcal{L}}
\newcommand{\tskip}{t_\mathrm{skip}}

\newcommand{\tauskip}{\tau_\mathrm{skip}}
\newcommand{\deltat}{\delta}
\newcommand{\Deltat}{\Delta}
\newcommand{\restrict}{\mathcal{R}}
\newcommand{\tup}{T_\mathrm{up}}
\newcommand{\dom}{\operatorname{dom}}
\newcommand{\cond}{\operatorname{cond}}
\newcommand{\avg}[1]{\langle #1 \rangle}
\newcommand{\norm}[1]{\| #1 \|}

\newcommand{\R}{\mathbb{R}}

\newcommand{\id}{I}

\newcommand{\rg}{\operatorname{rg}}

\renewcommand{\d}{\mathop{}\!\mathrm{d}}

\renewcommand{\phi}{\varphi}

\renewcommand{\epsilon}{\varepsilon}

\newtheorem{assumption}{Assumption}

\title{Implicit Methods for Equation-Free Analysis: Convergence
  Results and Analysis of Emergent Waves in Microscopic Traffic Models}

\author{Christian Marschler\thanks{Department of Mathematics,
    Technical University of Denmark, Matematiktorvet 303B, DK-2800
    Kgs. Lyngby, Denmark ({\tt c.marschler@mat.dtu.dk}).}
  \and Jan Sieber\thanks{College of Engineering, Mathematics and Physical Sciences,
    University of Exeter, North Park Road, Exeter (Devon) EX4 4QF
    ({\tt j.sieber@exeter.ac.uk}).}
  \and Rainer Berkemer\thanks{AKAD University of Applied Sciences, 
    Maybachstrasse 18-20, D-70469 Stuttgart, Germany 
    ({\tt rainer.berkemer@akad.de}).}
  \and Atsushi Kawamoto\thanks{Toyota Central R\&D Labs., Inc., Nagakute, Aichi 480-1192,
    Japan ({\tt atskwmt@mosk.tytlabs.co.jp}).}
  \and Jens Starke\thanks{Department of Mathematics,
    Technical University of Denmark, Matematiktorvet 303B, DK-2800 Kgs. Lyngby, Denmark
    ({\tt j.starke@mat.dtu.dk}).}
 }

\begin{document}

\maketitle

\newcommand{\slugmaster}{%
}

\begin{abstract}
  We introduce a general formulation for an implicit equation-free method
  in the setting of slow-fast systems.
  First, we give a rigorous convergence result for equation-free
  analysis showing that the implicitly defined coarse-level time stepper converges
  to the true dynamics on the slow manifold within an error that is
  exponentially small with respect to the small parameter measuring
  time scale separation. Second, we apply this result to the idealized
  traffic modeling problem of phantom jams generated by cars with
  uniform behavior on a circular road.  The traffic jams are waves
  that travel slowly against the direction of traffic. Equation-free
  analysis enables us to investigate the behavior of the microscopic
  traffic model on a macroscopic level. The standard deviation of
  cars' headways is chosen as the macroscopic measure of the
  underlying dynamics such that traveling wave solutions correspond to
  equilibria on the macroscopic level in the equation-free setup. The
  collapse of the traffic jam to the free flow then corresponds to a
  saddle-node bifurcation of this macroscopic equilibrium. We continue
  this bifurcation in two parameters using equation-free analysis.
\end{abstract}

\begin{keywords} 
equation-free methods, implicit methods, lifting, traffic modeling,
optimal velocity model, traveling waves, stability of traffic jams
\end{keywords}

\begin{AMS}
65P30, 
37M20, 
37Mxx, 
34E13 
\end{AMS}

\section{Introduction}
When one studies systems with many degrees of freedom, for example,
systems with a large number of particles or interacting agents, one is
often interested not so much in the trajectories at the microscopic
level (that is, of individual particles), but in the behavior on the
macroscopic scale (of the overall distribution of particles). The
classical example is the motion of molecules of a gas, resulting in
the laws of thermodynamics. In this classical case the macroscopic
description is derived in statistical mechanics from knowledge about
the microscopic behavior through time scale separation. 
 Other important examples are emerging patterns in physical,
  chemical, and biological systems, \eg Rayleigh-B\'enard convection
  rolls \cite{Rayleigh1916}, the Belousov-Zhabotinsky reaction
  \cite{Belousov1959,Zhabotinsky1964}, and stripes on zebra skin or
  patterns on butterfly wings \cite{Young1984}. A common approach in
the physics literature to deriving macroscopic descriptions are the
so-called adiabatic elimination or the slaving principle
\cite{Haken1983,Haken1983a}. These concepts are related to the
theorems  in the
mathematical literature about reductions to center manifolds or slow manifolds \cite{E2003,Kelley1967,Vanderbauwhede1989}.

For systems where no explicit macroscopic description can be derived
from microscopic models, Kevrekidis
and coworkers proposed that, if the number of particles is moderate,
then it is sometimes possible to skip the derivation of a macroscopic
description 
by performing the analysis of the dynamics in the macroscopic scale
directly. This approach relies on evaluating short bursts of
appropriately initialized simulations of the microscopic model (see,
for example,
\cite{kevrekidisgear2003,KevrekidisSamaey2009,Kevrekidis2010} for
recent reviews). It is called \emph{equation-free} because it assumes
that the macroscopic model exists but is not available as an explicit
formula. Equation-free methods are particularly appealing if either
explicit macroscopic descriptions are unavailable, or one wants to
study the underlying system near the boundary of validity of its
macroscopic description (for example, as one decreases the number of
particles, finite size effects may start to appear as small
corrections to the macroscopic model).  Equation-free analysis has
been applied for a large class of multiscale models that roughly fit
the description of \emph{singularly perturbed systems}
\cite{Fenichel1979} in a broad sense (see motivation in
\cite{Kevrekidis2010}), such as stochastic systems
\cite{Makeev2002,Siettos2003}, agent-based models
\cite{Cisternas2004,Corradi2012,0295-5075-82-3-38004}, molecular
dynamics \cite{chen_debendetti} or neural dynamics
\cite{springerlink:10.1007/s10827-006-3843-z,Reppas2010}, to perform
high-level tasks such as bifurcation analysis, optimization or control
design \cite{Elmegaard2009,AIC:AIC690490727}.

\begin{figure}[t]
  \centering
 \includegraphics{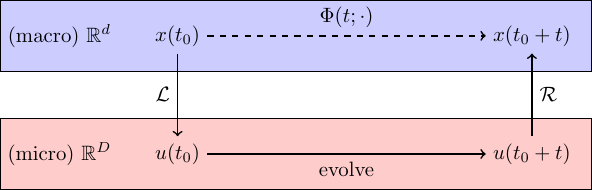}
  \caption{Sketch of the macroscopic time stepper $\Phi(t;\cdot)$. The macroscopic
    state $x(t_0)$ is mapped to a microscopic state $u(t_0)$ by using
    the lifting operator $\lift$. The available microscopic time
    stepper is used to evolve the system to the microscopic state
    $u(t_0 + t)$, which is mapped to a macroscopic state $x(t_0 + t)$
    using the restriction operator $\restrict$. This procedure
    constitutes the coarse-level time stepper $\Phi(t;\cdot)$.}
  \label{fig:sketch-eqfree}
\end{figure}

The basic building block of equation-free analysis is an approximate
coarse-level time stepper $\Phi(t;\cdot)$ for short times $t$
(compared to the slow time scale) in the phase space of macroscopic
variables (say, $\R^d$). This coarse-level time stepper is typically
composed of three steps: \emph{lift} (operator $\lift$), \emph{evolve},
and \emph{restrict} (operator $\restrict$), as shown in Figure
\ref{fig:sketch-eqfree}. To compute the map $\Phi(t;x)$ on a given
macroscopic state $x\in\R^d$, one has to apply a lifting operator
$\lift$ to map $x$ to a microscopic state $u\in\R^D$ (typically, $D\gg
d$); then one runs the microscopic simulation for the time $t$; and
finally one maps the end state of the microscopic simulation back into
$\R^d$ using a restriction operator $\restrict$.  A proof of any claim
that this would be a good approximation of the true dynamics of the
macroscopic variable $x$ for a given example will have to invoke the
following sequence of arguments.  Initially assume that the
microscopic system is a slow-fast system with a transversally stable
slow manifold, for which the macroscopic quantity $x$ is a
coordinate. The first question is then: does the approximate
coarse-level time stepper $\Phi$ converge to the true dynamics on the
slow manifold in the limit $\epsilon\to0$, where $\epsilon$ is the
parameter measuring the time scale separation?  In addition to the
case discussed here, equation-free analysis is also applied to
high-dimensional, stochastic (or chaotic) systems showing macroscopic
behavior because the dynamics of the microscopic degrees of freedom
averages out rapidly \cite{Barkley2006,Spohn1991,SanchezPalencia1983}.
In these cases another question must be addressed: in which sense is
the averaging process approximating a classical slow-fast system?

\subsection{An implicit coarse-level time stepper}
Before equation-free analysis can be performed, one must find the
restriction and lifting operators $\restrict$ and
$\lift$. Figure~\ref{fig:sketch-eqfree} suggests the relation
$\Phi(t;\cdot) = \restrict \circ \text{evolve} \circ \lift$. However,
this will not approximate the true macroscopic flow in general. Why?
Let us assume that the microscopic system is slow-fast and the
macroscopic system corresponds to the slow flow on the slow manifold
in the coordinate $x$. Then an arbitrary choice of $\lift$ and
$\restrict$ does not lead to a coarse time-stepper $\Phi$ which
approximates the slow flow in any way, even in the limit of infinite
time scale separation ($\epsilon\to0$). The source of the error is an
initialization of the microscopic system away from the slow manifold.
One relies on the separation of time scales in a so-called healing
step to reduce this error. However, in most reviews this healing is
applied inconsistently
\cite{kevrekidisgear2003,KevrekidisSamaey2009,Kevrekidis2010}. That
is, healing would not lead to $\Phi$ converging to the true slow flow
in the limit of infinite time scale separation, even in the ideal case
of a slow-fast system. A consistent way to perform healing are
so-called \emph{constrained-runs} corrections after lifting, developed
in \cite{Gear2005,Zagaris2009,Zagaris2012}. These papers developed
schemes of increasing complexity to compensate for this error source.

An alternative, explained in Section \ref{sec:implicit}, is to use an
implicitly defined coarse-level time stepper $\Phi$, where the slow
flow is not measured at predetermined points in space but rather at
healed points. In the special case of computation of equilibria, the
use of the implicit time stepper reduces to the formula introduced as
the ``third method'' by Vandekerckhove \emph{et al}
\cite{Vandekerckhove2011}.  In Section \ref{sec:thm}, we give a
detailed proof of the convergence of the implicitly defined
coarse-level time stepper $\Phi$ to the flow on the slow manifold,
answering the question of convergence for the implicit time
stepper. The approximation error of $\Phi$ (under some transversality
conditions) is exponentially small in the parameter $\epsilon$
measuring the time scale separation. Our theorem does not require
  that the time scale separation parameter $\epsilon$ approach zero,
  merely that it be sufficiently small. The precise statement is then
  that the error is of order $\exp(-K\tskip)$, where $K$ is
  the rate of attraction transversal to the slow manifold and $\tskip$
  is the healing time. In Section \ref{sec:comp} we discuss the
assumptions and consequences of the convergence theorem and compare it
to other results in the literature.

\subsection{Macroscopic behavior of a microscopic traffic model}
\label{sec:intro:traffic}

In Section \ref{sec:model} and Section \ref{sec:eqfree} we apply the
implicit coarse time stepper to a traffic modeling problem that fits
into the framework of equation-free analysis: a large number of cars
(the microscopic particles) on a circular road that interact with each
other, resulting in so-called phantom jams moving slowly along the
road against the direction of traffic, \ie forming a traveling wave at
the microscopic level.

The mathematical modeling and analysis of traffic flow dynamics has a
considerable history (see, \eg \cite{Helbing2001,Nagatani2002,Orosz2010}
for reviews). 
Macroscopic traffic models use partial differential equations, such as
Burger's equation \cite{Nagatani2002}, for modeling the flow.  They
model the density of cars as a continuous quantity to directly formulate
macroscopic equations for density and flux along the road. In
contrast, microscopic particle models (deterministic
\cite{Bando1995} or stochastic \cite{Kanai2009,Shigaki2011}) can be
used to describe the behavior of individual cars or drivers. An
advantage of microscopic models is that parameters can be assigned
directly to the individual drivers' behavior (for example,
aggressiveness, inertia, or reaction delay) such that these parameters'
influence and the trajectories of individual cars can be
investigated. Another use of microscopic models is to test the effects
of new devices for individual cars, for example, cruise control, on
the overall traffic prior to their implementation in real traffic. In
this paper we use the optimal velocity model \cite{Bando1995} as an
example of an underlying microscopic model. The optimal velocity model
results in a set of coupled ordinary differential equations, but
despite its simplicity it can reproduce the phenomenon of phantom
traffic jams. An advantage of choosing the optimal velocity model is
that we have guidance from the results of direct bifurcation analysis
of the full microscopic system when only a few cars are involved
\cite{Gasser2004,Orosz2005} as well as from perturbation analysis
based on the discrete modified Korteweg--de Vries equation
\cite{Gaididei2009}. Direct bifurcation analysis of the
microscopic system becomes infeasible when the number of cars gets
large. Furthermore, it is difficult to analyze macroscopic quantities
for which typically no equations are explicitly given such as the mean
and standard deviation of
headways or densities of cars. In Section \ref{sec:eqfree} we show how
this difficulty can be tackled by using equation-free methods for the
bifurcation analysis on a macroscopic level.

In Section \ref{sec:discussion} we summarize the
obtained results and give an overview of open problems.

\section{Nontechnical description of general equation-free analysis
  with implicit lifting}
\label{sec:implicit}
Equation-free analysis as described by \cite{KevrekidisSamaey2009} is
motivated by ideas from the analysis of slow-fast systems: one assumes
that on a long time scale the dynamics is determined by only a few
state variables and the other state variables are
\emph{slaved}. Mathematically this means that the flow of a
high-dimensional system under study converges rapidly onto a
low-dimensional manifold on which the system is governed by an
ordinary differential equation (ODE). In many practical applications
convergence is achieved only in the sense of statistical mechanics
(the effects of many particles averaging out; see
\cite{Barkley2006,Corradi2012}). We give our description
  and subsequent convergence proofs of equation-free analysis using
  the terminology of slow-fast systems with transversally stable slow
  manifolds following the notation of \cite{Fenichel1979}. 
  The
  traffic problem discussed in Section~\ref{sec:model} and
  \ref{sec:eqfree} does not require the notion of weak (averaged)
  convergence.

\subsection{The notion of a slow-fast system}
\label{sec:slowfastnotion}
Let
\begin{equation}\label{eq:f}
  \dot u=f_\epsilon(u)
\end{equation}
be a smooth dynamical system defined for $u\in\R^D$, where $f_\epsilon$ depends
smoothly on the parameter $\epsilon$. We assume that $\epsilon$ is a
\emph{singular perturbation parameter}. This means that the flow
$M_\epsilon$ generated by \eqref{eq:f},
\begin{displaymath}
M_\epsilon:\R\times \R^D\to \R^D\mbox{,} \qquad (t;u)\mapsto M_\epsilon(t;u)
\end{displaymath}
has a whole smooth $d$-dimensional submanifold ${\cal C}_0$ of
equilibria for $\epsilon=0$: if $u\in {\cal C}_0$, then $M_0(t;u)=u$
(and, thus, $f_0(u)=0$) for all $t$. The dimension $d$ is the number
of slow variables. In the notation of singular perturbation theory,
$t$ measures the time on the \emph{fast} time scale. We assume that
this manifold ${\cal C}_0$ is transversally uniformly exponentially
stable for $\epsilon=0$, which corresponds to the stable case of
Fenichel's geometric singular perturbation theory
\cite{Fenichel1979}. For this case we know that the flow
$M_\epsilon(t;\cdot)$ has a transversally stable invariant manifold
${\cal C}_\epsilon$ for small nonzero $\epsilon$, too. This manifold
${\cal C}_\epsilon$ is called the slow manifold, and the flow
$M_\epsilon$, restricted to ${\cal C}_\epsilon$, is called the slow
flow.  For the traffic problem the time scale separation is
  present as demonstrated numerically later in Section~\ref{sec:tsep}.

\subsection{Lifting, restriction, and time stepping}
The equation-free approach to coarse graining
\cite{KevrekidisSamaey2009} does not require direct access to the
right-hand side $f_\epsilon$ of the microscopic system \eqref{eq:f}
but merely the ability to evaluate $M_\epsilon(t;u)$ for finite
positive times $t$ (typically $t\ll 1/\epsilon$ in the fast time scale
$t$) and arbitrary $u$. It also relies on two smooth maps that have to
be chosen beforehand:
\begin{align*}
  \restrict&: \R^D\to \R^d &&\mbox{the \emph{restriction} operator,}\\
  \lift&: \R^d\to\R^D &&\mbox{the \emph{lifting} operator.}
\end{align*}
In the optimal velocity model discussed in Section \ref{sec:eqfree}, $\restrict$ is
chosen as a mapping from headway profiles to the standard deviation $\sigma$
and $\lift$ constructs a headway profile by using $\sigma$
(cf. \eqref{eq:restriction} and \eqref{eq:reslift}).
\par
\begin{figure}[t]
  \centering
  \subfigure[{implicit
    scheme}]{\label{sfig:sketch_impl}\includegraphics[scale = 0.72
    ]{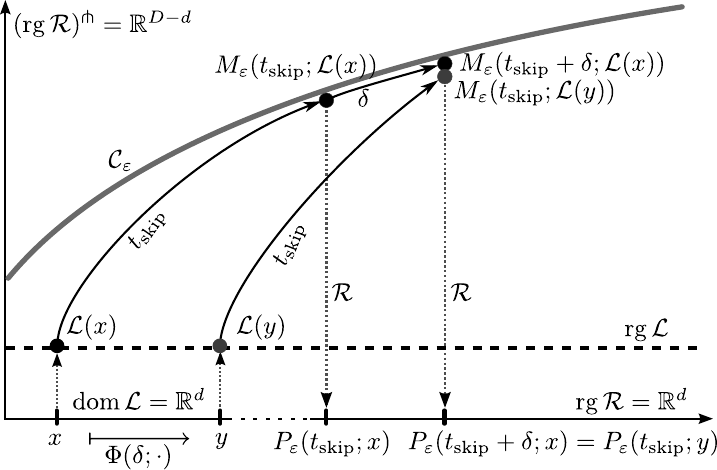}} 
\hfill 
\subfigure[{explicit
    scheme}]{\label{sfig:sketch_expl}\includegraphics[scale = 0.72
    ]{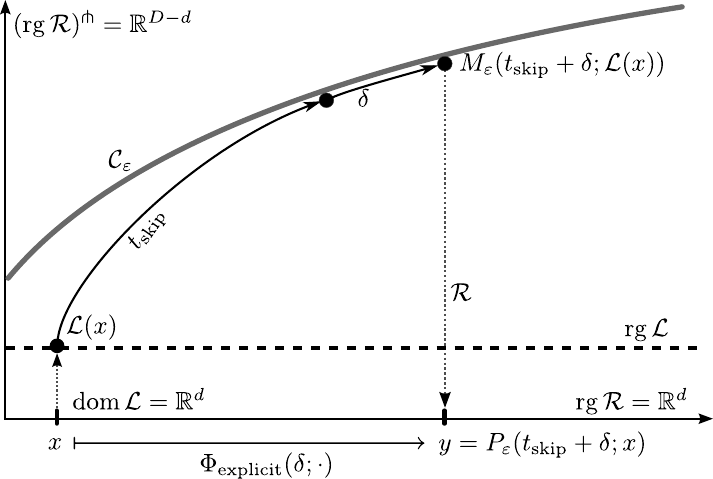}}

  \caption{\textup{(}a\textup{)} Sketch showing a typical geometry of
    the implicit scheme in a slow-fast system with a
    slow manifold ${\cal C}_\epsilon$ and an arbitrary lifting $\lift$
    and restriction $\restrict$. The healing
    $M_\epsilon(\tskip;\cdot)$ is applied to all points in the domain
    of $\lift$. Note that $\dom\lift$ and $\rg\restrict$ can be
    different, but must have the same dimension. $(\rg
    \restrict)^\pitchfork$ refers to an arbitrary transversal
    complement of $\rg\restrict$. \textup{(}b\textup{)} The explicit
    scheme is shown for comparison.}
  \label{fig:sketch-overview}
\end{figure}

 
The basic idea underlying \cite{KevrekidisSamaey2009} is that one can
analyze the dynamics of \eqref{eq:f} on the slow manifold ${\cal
  C}_\epsilon$ by studying a map in the space of restricted variables
$x$ in the domain of $\lift$ (called $\dom \lift\subset\R^d$) of the
form (cf. Figure \ref{fig:sketch-eqfree})
\begin{center}
  \emph{Lift}\quad$\to$\quad\emph{Evolve}\quad$\to$\quad\emph{Restrict},
\end{center}
or, to be  precise, the map
\begin{equation}
  \label{eq:eqfree:Peps}
  P_\epsilon(t;\cdot):x\mapsto
  \restrict(M_\epsilon(t;\lift(x)))=[\restrict\circ M_\epsilon(t;\cdot)\circ \lift](x)
\end{equation}
for selected times $t\ll1/\epsilon$.  The central question is: how can
one compose a macroscopic time stepper, that is, an approximate
time-$\deltat$ map $\Phi(\deltat,\cdot):\R^d\to \R^d$, using
coordinates in the domain of $\lift$ for the flow $M_\epsilon$
restricted to ${\cal C}_\epsilon$? One important observation is that
this map $\Phi$ must be defined
implicitly. Figure~\ref{fig:sketch-overview} shows how one can define
a good approximate time-$\deltat$ map $\Phi(\deltat;\cdot)$. It
contains an additional parameter $\tskip$, called the \emph{healing
  time} in \cite{KevrekidisSamaey2009}.
This healing time must be applied to both the argument $x$ and the
result $y$ of $\Phi$. Thus, $\Phi(\deltat;\cdot)$ is given implicitly
by solving
\begin{equation}
    \begin{aligned}
      P_\epsilon(\tskip;y)=&P_\epsilon(\tskip+\deltat;x)\mbox{,}
      &&\mbox{that is,}\\
      \restrict(M_\epsilon(\tskip;\lift(y)))
      =&\restrict(M_\epsilon(\tskip+\deltat;\lift(x)))\mbox{} &&\mbox{}
    \end{aligned}
    \label{eq:phidef}
\end{equation} 
for $y$, and setting $\Phi(\deltat;x):=y$.
Under some genericity conditions on $\restrict$, $\lift$, and
$M_\epsilon$ the order of approximation for $\Phi$ is exponentially
accurate for increasing $\tskip$ if we assume that $\epsilon\tskip$
and $\epsilon(\tskip+\deltat)$ are bounded:
\begin{equation}
  \label{eq:accuracy}
  \|\Phi(\deltat;x)-\Phi_*(\deltat;x)\|\leq 
  C\exp(-K\,\tskip)\mbox{.}
\end{equation}
In this estimate $K>0$ and $C>0$ are constants that depend only on a
uniform upper bound $\tup$ for $\epsilon\tskip$ and
$\epsilon(\tskip+\deltat)$.
The flow $\Phi_*$ is the exact flow $M_\epsilon$, restricted to the
slow manifold ${\cal C}_\epsilon$, in a suitable coordinate
representation in $\dom\lift$. The same estimate holds also for the
derivatives of $\Phi$ with respect to the initial value up to a fixed
order (with more restrictive conditions on $\epsilon$). So,
\begin{displaymath}
  \|\partial_2^j\Phi(\deltat;x)-\partial_2^j\Phi_*(\deltat;x)\|\leq 
C\exp(-K\tskip)
\end{displaymath}
(possibly with other constants $C$) for derivative orders $j$ less
than a given $k$ (the subscript of $\partial_i^j$ refers to the
argument of $\Phi$ with respect to which the $j$th derivative is
taken). The degree of achievable differentiability is determined by
the time scale separation: the smaller $\epsilon$ is, the smoother the
slow manifold ${\cal C}_\epsilon$ is, and, thus, the higher we can choose
the maximal derivative order $k$.

Based on the implicitly defined approximate flow map $\Phi$, one can now perform
higher-level tasks in equation-free analysis.

\subsection{Bifurcation analysis of macroscopic equilibria}
Bifurcation analysis for equilibria boils down to finding fixed points
and their stability and bifurcations for $\Phi(\deltat;\cdot)$ with
some small, arbitrary $\deltat$ (that is,
$\deltat\ll1/\epsilon$ in our notation). In terms of $\restrict$ and
$\lift$, the equation $\Phi(\deltat;x_0)=x_0$, defining the
equilibrium $x_0$, reads (cf. Figure \ref{fig:sketch-overview})
\begin{equation}
  \label{eq:eqfree:simple}
  \restrict(M_\epsilon(\tskip+\deltat;\lift(x_0)))=\restrict(M_\epsilon(\tskip;\lift(x_0)))\mbox{.}
\end{equation}
This equation has been proposed and studied already in
\cite{Vandekerckhove2011}.
In applications, \eqref{eq:eqfree:simple} is solved using a
Newton iteration (cf. \eqref{eq:newton} in the optimal velocity model).
  Since the time stepper is defined
implicitly, one finds the stability
and bifurcations of an equilibrium $x_0$ by studying the generalized
eigenvalue problem
\begin{equation}\label{eq:matrixpair}
  \left[\frac{\partial}{\partial x}  
    \left[\restrict(M_\epsilon(\tskip+\deltat;\lift(x)))\right]\Bigl\vert_{x=x_0}\right] x=
    \lambda\left[\frac{\partial}{\partial x} \left[\restrict(M_\epsilon(\tskip;\lift(x)))
    \right]\Bigl\vert_{x=x_0}\right]x
  \mbox{.}
\end{equation}
This eigenvalue problem will give the eigenvalues of the
implicitly-known flow $\Phi(\delta;\cdot)$, linearized with respect to
its second argument $x$ in the equilibrium $x_0$ such that
bifurcations occur when $\lambda$ is on the unit circle.
\subsection{Projective integration}
\label{sec:proj:int}
In projective integration one approximates the ODE for the flow on the
slow manifold ${\cal C}_\epsilon$ in the coordinate $x\in\R^d$. The
ODE for the true flow $\Phi_*$ on the slow manifold is an implicit ODE
with the solution $x(t)$, which will be derived in detail in Section
\ref{sec:thm}. Its approximation based on $\Phi$ is
\begin{equation}
  \label{eq:implicitode}
  \frac{\d}{\d t}\restrict(M_\epsilon(\tskip;\lift(x)))=
  \frac{\partial}{\partial\deltat}
  \restrict(M_\epsilon(\tskip+\deltat;\lift(x)))\Bigl\vert_{\deltat=0}\mbox{.}
\end{equation}
For fixed $\tskip$ the left-hand side is a function of $x\in\R^d$ such
that the time-derivative of this function defines (implicitly) the
time-derivative of $x$. The term inside the partial derivative on the
right-hand side is a function of two arguments, $\deltat$ and $x$, for
which one takes the partial derivative with respect to its first
argument $\deltat$ in $\deltat=0$, making also the right-hand side a
function of $x$ only. Consequently, every integration scheme becomes
implicit. For example, if one wants to perform an explicit Euler step
of stepsize $\Delta t$ starting from $x_j$ at time
$t_j$, this becomes an implicit scheme (defining $x_{j+1}$ as the new
value at time $t_{j+1} = t_j + \Delta t$):
\begin{equation}
  \label{eq:expleuler}
  \frac{1}{\Delta t}\left[P_\epsilon(\tskip;x_{j+1})-P_\epsilon(\tskip;x_j)\right]=
  \frac{1}{\deltat}\left[P_\epsilon(\tskip+\deltat;x_j)-
    P_\epsilon(\tskip;x_j)\right]\mbox{,}
\end{equation}
or, in terms of restricting and lifting,
\begin{multline*}
  \restrict(M_\epsilon(\tskip;\lift(x_{j+1})))-
  \restrict(M_\epsilon(\tskip;\lift(x_j)))\\
  =
  \frac{\Delta t}{\deltat}\left[\restrict(M_\epsilon(\tskip+\deltat;\lift(x_j)))-
    \restrict(M_\epsilon(\tskip;\lift(x_j)))\right]
\end{multline*}
Projective integration becomes attractive if either one can choose
$\Delta t$ much larger than $\tskip$ and
$\deltat$, or one can set $\Delta t$ negative, enabling
integration backward in time on the slow manifold
(cf. \eqref{eq:bweuler} and Figure \ref{fig:profiles}), even though
the original system is very stiff in $\R^D$ forward in time (and
thus, strongly expanding backward in time). For positive $\Delta
  t$ the restriction on the size of $\Delta t$ is given by standard
  consistency and stability requirements of the coarse-grained
  integration method restricted to the slow flow (in general the
  restriction is $\epsilon\Delta t\ll1$, which makes the maximal stepsize
  independent of the time-scale separation). 
Note that during
computation of residuals and Jacobian matrices one can evaluate
$P_\epsilon(\tskip;x)$ as a by-product of the evaluation of
$P_\epsilon(\tskip+\deltat;x)$, assuming that the restriction
$\restrict$ is of comparatively low computational cost.

\subsection{Matching the restriction}
Sometimes it is of interest to find a microscopic state $u\in\R^D$ on
the slow manifold ${\cal C}_\epsilon$ that has a particular $x\in\R^d$
as its restriction ($\restrict(u)=x$); see
\cite{Gear2005,Zagaris2009,Zagaris2012}. This state $u$ is defined
implicitly and can be found by solving the $d$-dimensional nonlinear
equation
\begin{equation}
  \label{eq:rmatch}
  \restrict(M_\epsilon(\tskip;\lift(\tilde x)))=x
\end{equation}
for $\tilde x$, and then setting $u=M_\epsilon(\tskip;\lift(\tilde
x))$. This solution $u$ is close to the true slow manifold ${\cal
  C}_\epsilon$ with an error of order 
$\exp(-K\tskip)$, where the decay rate $K>0$ and
the possible constant in front of the exponential are independent of
$\epsilon$ and $\tskip$. This implies that, if we choose
$\tskip=O(\epsilon^{-1})$ with $p\in(0,1)$, the distance of $u$ to
${\cal C}_\epsilon$ is small beyond all orders of $\epsilon$ (see
Section~\ref{sec:thm} for the precise conditions). Equation
\eqref{eq:rmatch} was also proposed and studied in
\cite{Vandekerckhove2011} (called \textsc{InitMan} in
\cite{Vandekerckhove2011}), although without the general error estimate.

\section{Convergence of equation-free analysis}
\label{sec:thm}
This section gives a detailed discussion of the convergence results of
the methods sketched in Section \ref{sec:implicit}. Sections
\ref{sec:model} and \ref{sec:eqfree} study the optimal
  velocity model for traffic flow as an application of implicit
equation-free analysis.

We formulate all assumptions on $\restrict$, $\lift$, and $M_\epsilon$
for the singular perturbation parameter $\epsilon$ at $\epsilon=0$,
even though it is typically difficult to vary $\epsilon$ in complex
model simulations. However, stating the conditions at $\epsilon=0$
ensures that they are uniformly satisfied for all sufficiently small
$\epsilon$, which is the range of parameters for which the statements
of this section are valid (cf. \cite{Fenichel1979}). Throughout this
section various constants will appear in front of exponentially
growing or decaying quantities. As the concrete values of these
constants do not play a role, we will use the same variable name $C$ on
all occasions without meaning them to be the same. We will state which
quantities the constant $C$ depends on whenever we use exponential
estimates.

The notation $\partial_k^j$ refers to the $j$th derivative with
respect to the $k$th argument. For example, $\partial_2^jM_\epsilon$
refers to the $j$th-order partial derivative of the flow $M_\epsilon$
with respect to its second argument (the starting point), and the
zeroth derivative refers to the value of flow $M_\epsilon(t;\cdot)$
itself.
\begin{figure}[ht]
  \centering
  \includegraphics[width=0.9\textwidth]{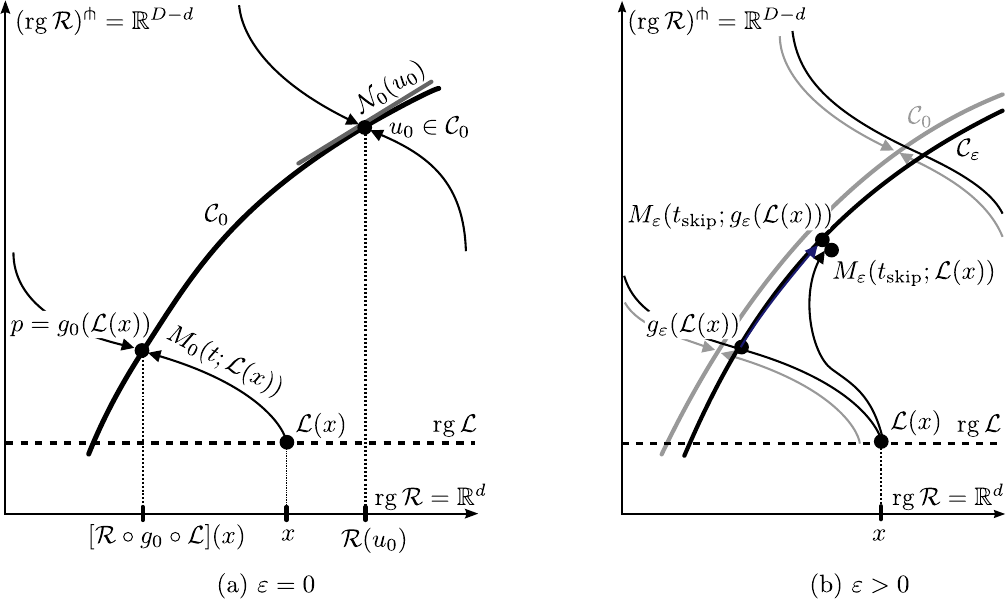}
  \caption{Sketch of geometrical interpretation of transversality
    assumptions. Note that $(\rg \restrict)^\pitchfork$ refers to an arbitrary
    transversal complement of $\rg\restrict$. Panel
    \textup{(}a\textup{)}  shows the geometry
    at $\epsilon=0$: The trajectory
    starting at
    $\lift(x)$ must converge to ${\cal C}_0$, and its limit is called
    $g_0(\lift(x))$. The overall map $\restrict\circ g_0\circ \lift$
    must be a local diffeomorphism. This entails that
    the Jacobian $\partial\restrict$ must have full rank on the tangent space
    ${\cal N}_0(u_0)$ in any $u_0\in{\cal C}_0$ \textup{(}also shown
    in \textup{(}a\textup{))}, and that $\rg\lift$ intersects each
    fiber \textup{(}the set of points $u$ converging to the same
    $u_0\in{\cal C}_0$\textup{)} transversally. Shown in
    panel \textup{(}b\textup{)}: $g_\epsilon$ and ${\cal C}_\epsilon$
    are $O(\epsilon)$ perturbations of $g_0$ and ${\cal C}_0$, and
    $M_\epsilon(\tskip;\lift(x))-M_\epsilon(\tskip;g_\epsilon(\lift(x)))$
    are $\exp(-K \tskip)$ close for $\tskip>0$.}
  \label{fig:sketch-geom}
\end{figure}
\subsection{Existence of transversally stable slow manifold}
As introduced in Section~\ref{sec:slowfastnotion}, the
  microscopic flow $M_\epsilon(t;u_0)$ is the solution of
  \begin{equation}
    \label{eq:frep}
    \dot u=f_\epsilon(u)\mbox{,\quad $u(t)\in\R^D$,}
  \end{equation}
  starting from initial condition $u_0\in\R^D$, which for $\epsilon=0$ has a
  $d$-dimensional manifold of equilibria ${\cal C}_0$. That is,
  $f_0(u)=0$ if $u\in{\cal C}_0$.  In order to avoid the discussion of
what happens when the flow $M_\epsilon$ reaches certain boundaries or
  becomes large while following the slow dynamics, we assume that the
manifold ${\cal C}_0$ of equilibria of $M_0$ is compact.  Our first
assumption guarantees transversal stability of ${\cal C}_0$.
\begin{assumption}[Separation of time scales and transversal
  stability]\label{ass:timescales}
  \quad There exists a constant $K_0>0$ such that for all points $u\in{\cal
    C}_0$ the Jacobian $\partial f_0(u)$ has $D-d$ eigenvalues with
  real part less than $-K_0$.
\end{assumption}

This implies that the flow $M_0$
approaches the slow manifold ${\cal C}_0$
with a rate faster than $K_0$ from all initial conditions $u$ in some
neighborhood of ${\cal C}_0$. That is, for every $u$ in an
appropriate open neighborhood ${\cal U}$ of the slow manifold
${\cal C}_0$ there exists a point $p\in{\cal C}_0$ such that
\begin{displaymath}
  \lim_{t\to\infty}M_0(t;u)=p
\end{displaymath}
(note that for $\epsilon=0$ all points on the slow manifold ${\cal
  C}_0$ are equilibria), and the distance can be bounded via
\begin{displaymath}
  \|M_0(t;u)-p\|\leq C\exp(-K_0t)\|u-p\|\mbox{,\quad}
  \|\partial_2^jM_0(t;u)\|\leq C\exp(-K_0t)
\end{displaymath}
for all $t\geq0$ and $j\geq1$, where the constant $C$ depends only
on the derivative order $j$.

Since the slow manifold ${\cal C}_0$ is compact, one can choose a
uniform constant $C$ for all $u$ in the neighborhood ${\cal U}$. The
above assumption implies the existence of a smooth map (called the
\emph{stable fiber projection}),
\begin{equation}
  g_0:{\cal U}\to{\cal C}_0\mbox{,\quad defined by\quad} 
  g_0(u):=p\mbox{,}\label{eq:gdef}
\end{equation}
assigning to each $u$ its limit $p\in{\cal C}_0$ under the flow $M_0$ (see
Figure~\ref{fig:sketch-geom}(a)).

We recall now two central persistence results of classical singular
perturbation theory \cite{Fenichel1979}. First, the slow manifold
${\cal C}_0$ persists for sufficiently small $\epsilon$, deforming to
a smooth nearby manifold ${\cal C}_\epsilon$ (as shown in
Figure~\ref{fig:sketch-geom}(b)). This manifold ${\cal C}_\epsilon$ is
also compact. Restricted to ${\cal C}_\epsilon$, the flow $M_\epsilon$
is governed by a smooth ODE (the \emph{slow} flow) with a right-hand
side for which all derivatives up to a given order $k$ are
proportional to $\epsilon$ (larger $k$ requires smaller $\epsilon$):
\begin{equation}\label{eq:msmall}
\left\|f_\epsilon(u)\right\|\leq \epsilon\mbox{,\qquad}
  \left\|\partial^jf_\epsilon(u)[v_1,\ldots,v_j]\right\|\leq \epsilon\|v_1\|\cdot\ldots\cdot\|v_j\|
\end{equation}
for all $j=\{1,\ldots,k\}$, $u\in{\cal C}_\epsilon$ and
$v_1,\ldots,v_j\in{\cal N}_\epsilon(u)$. (Here ${\cal N}_\epsilon(u)$ is the
tangent space of ${\cal C}_\epsilon$; for $\epsilon=0$ it is the null
space of the linearization of $f_0$ in 
$u$ on the slow manifold ${\cal C}_0$.) Note that typically one has only
$\|\partial^jf_\epsilon(u)\vert_{{\cal C}_\epsilon}\|\leq C\epsilon$,
but we can set the constant $C$ equal to unity without loss of
generality by rescaling time or redefining the parameter $\epsilon$.
Thus, the flow $M_\epsilon(t;\cdot)$ is a global diffeomorphism on the
slow manifold ${\cal C}_\epsilon$ which has growth bounds of order
$\epsilon$ forward and backward in time:
\begin{align}
  \label{eq:mcebound}
    \|\partial_2^jM_\epsilon(t;\cdot)\vert_{{\cal C}_\epsilon}\|
    &\leq C\exp(\epsilon|t|)\mbox{,}&
  \|\partial_2^jM^{-1}_\epsilon(t;\cdot)\vert_{{\cal C}_\epsilon}\|
  &\leq C\exp(\epsilon|t|)\mbox{,}
\end{align}
for some constant $C$ independent of $t$ and $\epsilon$ and all
derivative orders $j$ up to a fixed order $k$. Note that
$M_\epsilon^{-1}(t;\cdot)=M_\epsilon(-t;\cdot)$ exists for all
times $t$ as long as one restricts the flow $M_\epsilon$ to the slow
manifold ${\cal C}_\epsilon$.

Second, the stable fiber projection map $g_0$ persists for small
$\epsilon$, getting perturbed smoothly to a map $g_\epsilon$, defined
for each $u$ in the neighborhood ${\cal U}$ of the slow manifold
${\cal C}_0$ (and its perturbation ${\cal C}_\epsilon$). The map
$g_\epsilon$ picks for every point $u\in{\cal U}$ the unique point
$g_\epsilon(u)$ inside the slow manifold ${\cal C}_\epsilon$ such that
the trajectories starting from $u$ and $g_\epsilon(u)$ converge to
each other forward in time with an exponential rate $K$ of order $1$
(that is, $K$ is uniformly positive for all sufficiently small
$\epsilon$ and all $u\in{\cal U}$):
\begin{equation}
  \label{eq:gepsconv}
  \|\partial_2^jM_\epsilon(t;u)-\partial_2^jM_\epsilon(t;g_\epsilon(u))\|\leq 
  C\exp(-Kt)\|u-g_\epsilon(u)\|
\end{equation}
for all $t\geq0$, $u\in{\cal U}$, and $0\leq j\leq k$, where the
constant $C$ is uniform for $u\in{\cal U}$.  In general, the decay
rate $K$ has to be smaller than the rate $K_0$ asserted to
exist in Assumption~\ref{ass:timescales} for $\epsilon=0$.  More
precisely, for every rate $K<K_0$ there exists a range
$(0,\epsilon_0)$ of $\epsilon$ for which \eqref{eq:gepsconv}
holds. Choosing $\epsilon_0$ smaller permits one to choose $K$ closer
to $K_0$.  The stable fiber projection map $g_\epsilon$ is an
order-$\epsilon$ perturbation of $g_0$:
\begin{equation}
  \label{eq:gepspert}
  \|\partial^j g_\epsilon(u)-\partial^jg_0(u)\|\leq C\epsilon
\end{equation}
for all $j=\{0,\ldots,k\}$ and a constant $C$ that is uniform for all
$u\in{\cal U}$.  The black curves transversal to ${\cal C}_\epsilon$
in Figure~\ref{fig:sketch-geom}(b) illustrate the fibers, that is,
which points of ${\cal U}$ get mapped onto the same point in ${\cal
  C}_\epsilon$ under $g_\epsilon$. Note that the fibers are not
trajectories for $\epsilon>0$; rather they are $(D-d)$-dimensional manifolds.

\subsection{Transversality conditions on restriction and lifting}
One assumption on the restriction $\restrict$ and the
lifting $\lift$ is that they are both smooth maps.

Furthermore, we assume that the lifting operator $\lift$ maps some
bounded open set $\dom\lift\subset\R^d$ into the basin of attraction
${\cal U}$ of ${\cal C}_0$ for $\epsilon=0$. We will make all
convergence statements in this section for $x\in\dom\lift$.

We formulate the transversality conditions on $\restrict$ and $\lift$
with the help of the tangent space ${\cal N}_0(u)$ to the slow
manifold ${\cal C}_0$ in a point $u_0\in{\cal C}_0$, which is given as
\begin{equation} \label{eq:tangentdef}
  {\cal N}_0(u_0)=\ker \partial f_0(u_0)\mbox{.}
\end{equation}
Remember that the
stable fiber projection $g_0$ maps all $u\in{\cal U}$ onto the slow
manifold ${\cal C}_0$. The tangent
space ${\cal N}_\epsilon(u)$ to the perturbed slow manifold ${\cal
  C}_\epsilon$ in a point $u\in{\cal C}_\epsilon$ is a perturbation of
${\cal N}_0(u)$ of order $\epsilon$.
\begin{assumption}[Transversality of $\restrict$ and
  $\lift$]\label{ass:transversality}
  \begin{enumerate}
  \item\label{ass:ltrans} The map $g_0\circ\lift$ is a local
    diffeomorphism between $\dom\lift\subset\R^d$ and the slow
    manifold ${\cal C}_0$ for every $x\in\dom\lift$.

    Equivalently, the composition of the linearizations $\partial
    g_0(\lift(x))\in\R^{D\times D}$ and
    $\partial\lift(x)\in\R^{D\times d}$ has full rank for all
    $x\in\dom\lift\subset\R^d$.
  \item\label{ass:rtrans} The map $\restrict:{\cal U}\to\R^d$,
    restricted to the slow manifold ${\cal C}_0$, is a local
    diffeomorphism between ${\cal C}_0$ and $\R^d$ for every $u$ in
    some relatively open subset $\dom\restrict\cap{\cal C}_0$.

    Equivalently, the dimension of the space $\partial
    \restrict(u){\cal N}_0(u)$ equals $d$ for every
    $u\in\dom\restrict\cap {\cal C}_0$.
  \item\label{ass:rdef} The set $\dom\restrict\cap{\cal C}_0$ contains
    $g_0(\lift(\dom\lift))$ as a subset, and the boundary of
    $\dom\restrict\cap{\cal C}_0$ has a positive distance from the
    boundary of $g_0(\lift(\dom\lift))$.
  \end{enumerate}
\end{assumption}
Note that points \ref{ass:ltrans} and \ref{ass:rtrans} of
Assumption~\ref{ass:transversality} are generically satisfied in a
given $x\in\R^d$ or $u_0\in{\cal C}_0$. 
By convention we keep $\dom\lift$ and
$\dom\restrict$ such that the transversality conditions are uniformly
satisfied in $\dom\lift$ and $\dom\restrict$.  The assumption that
$\dom\restrict$ (the region where $\restrict$ satisfies 
Assumption \ref{ass:rtrans}) contains the set $g_0(\lift(\dom\lift))$
guarantees that the map $x\mapsto \restrict(g_0(\lift(x)))$ is locally
invertible for all $x\in \dom\lift$ and that its linearization is
uniformly regular in $\dom\lift$. All points of
Assumption~\ref{ass:transversality} and the invertibility of the slow
flow, $M_\epsilon(t;\cdot)$ restricted to the slow
manifold ${\cal C}_\epsilon$, can be combined to ensure that the map
\begin{equation}\label{eq:adef}
  \R^d\supseteq \dom\lift \ni x\mapsto 
  \restrict(M_\epsilon(t;g_\epsilon(\lift(x))))\in\R^d
\end{equation}
is locally invertible for all $\epsilon\in[0,\epsilon_0)$ and
for all times $t$ satisfying
  \begin{equation}
    |t|\leq \tup/\epsilon\label{eq:tup}
\end{equation}
for some constant $\tup$. The constant $\tup$ is independent of
$\epsilon$, $t$, and $x\in\dom\lift$. It is determined by the distance
between the boundaries of $\dom\restrict$ and
$g_\epsilon(\lift(\dom\lift))$. This distance is positive because of
point~\ref{ass:rdef} in 
Assumption~\ref{ass:transversality} and the fact that $g_\epsilon$ is
a small perturbation of $g_0$. Then the time it takes a trajectory on
${\cal C}_\epsilon$ to reach the boundary of $\dom\restrict$ starting
from $g_\epsilon(\lift(\dom\lift))$ is of order $1/\epsilon$ such that
we can introduce the constant $\tup$. All components of the map
\eqref{eq:adef} are locally invertible:
$g_\epsilon\circ\lift:\dom\lift\to{\cal C}_\epsilon$ by
Point~\ref{ass:ltrans} of Assumption~\ref{ass:transversality}
(transversality of $\lift$); $M_\epsilon(t;\cdot)$ is a diffeomorphism
on ${\cal C}_\epsilon$; and $\restrict$, restricted to ${\cal C}_0$
(and, hence, to ${\cal C}_\epsilon$), is also locally invertible due to
Point~\ref{ass:rtrans} of Assumption~\ref{ass:transversality}. For
$\epsilon=0$ the map \eqref{eq:adef} is independent of $t$. Moreover,
the norm of the derivative of the map \eqref{eq:adef} and its inverse
are also uniformly bounded if $|\epsilon t|\leq \tup$

\subsection{Map of exact flow $M_\epsilon$ into $\R^d$}
Next, we give a coordinate system and a constructive procedure that
maps the flow $M_\epsilon$, restricted to the slow manifold ${\cal
  C}_\epsilon$, back to $\R^d$. This kind of map is called a
``lifting'' of the flow $M_\epsilon$ on ${\cal C}_\epsilon$ to its
cover $\R^d$ in, \eg \cite{Fenichel1979}, but we do not use this term
here to avoid confusion with the lifting operation $\lift$, used in an
equation-free context (cf. for example
\cite{KevrekidisSamaey2009}). For any fixed $\tskip$ the following map
$X_\epsilon:\dom\lift\to {\cal C}_\epsilon$ introduces
  coordinates of (part of) ${\cal C}_\epsilon$ in $\dom\lift$:
\begin{displaymath}
  X_\epsilon(x)=M_\epsilon(\tskip;g_\epsilon(\lift(x)))\mbox{.}
\end{displaymath}
This map is locally invertible because $g_0\circ \lift$ is a local
diffeomorphism between $\dom\lift$ and ${\cal C}_0$ (and, hence,
$g_\epsilon\circ\lift$ is a diffeomorphism between $\dom\lift$ and
${\cal C}_\epsilon$ for small $\epsilon$), and
$M_\epsilon(\tskip;\cdot)$ is a global diffeomorphism on ${\cal
  C}_\epsilon$ (see \eqref{eq:mcebound}). Moreover, if
  $M_\epsilon(\tskip;g_\epsilon(\lift(x)))$ is in the interior of the
  domain of $\restrict$, then one can find, for a given
$u=X_\epsilon(x)\in{\cal C}_\epsilon$, a preimage $\tilde x\approx x$ of any
point $\tilde u\in{\cal C}_\epsilon$ close to $u$ by solving
\begin{equation}\label{eq:mapback}
  \restrict(X_\epsilon(\tilde x))=\restrict(\tilde u)
\end{equation}
for $\tilde x$.  This follows from Assumption~\ref{ass:transversality}
(transversality for $\restrict$). In particular,
point~\ref{ass:rdef} of Assumption~\ref{ass:transversality} gives the
bound on the range of $\tskip$ for which the linearization of
\eqref{eq:mapback} is regular: the trajectory $t\mapsto
M_\epsilon(t;g_\epsilon(\lift(x)))$ should not leave $\dom\restrict$
for $t\in[0,\tskip]$, which is guaranteed for
$\tskip<\tup/\epsilon$. By requiring $\tilde x\approx x$, the
preimage $\tilde x$ of $\tilde u$, defined by \eqref{eq:mapback},
becomes unique.

Let $x(\deltat)\in \dom\lift\subset\R^d$ be a trajectory of the flow
  $M_\epsilon$ on ${\cal C}_\epsilon$ in the coordinates defined by
  $X_\epsilon$. By definition, $x$ satisfies
  $X_\epsilon(x(\delta))=M_\epsilon(\delta;X_\epsilon(x(0)))$. As long
  as $X_\epsilon(x(\delta))$ is in the domain of $\restrict$, we can apply
  $\restrict$ to this identity to obtain
  \begin{equation}
    \begin{aligned}
      \restrict(X_\epsilon(x(\delta)))&=
      \restrict(M_\epsilon(\delta;X_\epsilon(x(0))))\mbox{,} \qquad
      \qquad \mbox{that is,}\\
      \restrict(M_\epsilon(\tskip;g_\epsilon(\lift(x(\delta)))))&=
      \restrict(M_\epsilon(\tskip+\delta;g_\epsilon(\lift(x(0)))))\mbox{}
    \end{aligned}
\label{eq:phistardef0}
\end{equation}
(inserting the definition of $X_\epsilon$). Hence, the flow
$M_\epsilon$ on ${\cal C}_\epsilon$, written in the coordinates
  $x\in\dom\lift$, 
satisfies the implicit ODE
\begin{equation}
  \label{eq:exactode}
  \frac{\d}{\d t}\restrict(M_\epsilon(\tskip;g_\epsilon(\lift(x))))=
  \frac{\partial}{\partial\deltat}
  \restrict(M_\epsilon(\tskip+\deltat;g_\epsilon(\lift(x))))\Bigl\vert_{\deltat=0}\mbox{}
\end{equation}
as long as $\epsilon\tskip<\tup$ and $\epsilon(\tskip+\deltat)<\tup$
such that the resulting trajectory $x(\deltat)$ stays in $\dom\lift$
and
$X_\epsilon(x(\deltat))=M_\epsilon(\tskip;g_\epsilon(\lift(x(\deltat))))$
stays in $\dom\restrict$. For different values of $\tskip$ we get
different coordinate representations of the same flow, all related to
the representation with $\tskip=0$ via the global diffeomorphism
$M_\epsilon(\tskip;\cdot)$ on ${\cal C}_\epsilon$, which is a
near-identity transformation if $\tskip\ll1/\epsilon$ (see
\eqref{eq:mcebound}).

Let us denote the flow corresponding to the trajectory $x(\delta)$ in
\eqref{eq:phistardef0} as
$\Phi_*(\deltat;\cdot):\dom\lift\to\dom\lift$. The flow $\Phi_*$
is generated by the ODE \eqref{eq:exactode}. If
$\epsilon\tskip\leq\tup$ and $\deltat\ll1/\epsilon$, this flow map
$\Phi_*(\deltat;\cdot)$ is defined implicitly by solving the following
system for $y_*$,
\begin{equation}
  \restrict(M_\epsilon(\tskip;g_\epsilon(\lift(y_*))))=
      \restrict(M_\epsilon(\tskip+\deltat;g_\epsilon(\lift(x))))\mbox{,}
\label{eq:phistardef}
\end{equation}
and setting $\Phi_*(\deltat;x):=y_*$.  The local invertibility of
$X_\epsilon$ guarantees that there is a solution $y_*$ close to $x$
and that the solution $y_*$ is unique in the vicinity of $x$.  For
larger $\deltat$, one breaks down the flow into smaller time steps
such that one can apply the local solvability at every step:
\begin{equation}\label{eq:splittimestep}
  \Phi_*(\deltat;x)=\Phi_*(\deltat/m;\cdot)^m [x]
\end{equation}
for sufficiently large integer $m$. This construction of $\Phi_*$
achieves a representation of the exact flow $M_\epsilon$ restricted to
${\cal C}_\epsilon$ that is globally unique on $\dom\lift$ for all
$\deltat$ with $\epsilon(\tskip+\deltat)\leq\tup$.

\subsection{Approximate flow map and its convergence}
We now define the approximate flow map $y=\Phi(\deltat;x)$. Its
definition is similar to \eqref{eq:phistardef}, in particular, it is
also implicit. To highlight where the difference between $y$ and $y_*$
comes from, we put the defining equation for $y_*=\Phi_*(\deltat;x)$
directly below the implicit definition of $y$:
\begin{equation}\label{eq:exactperturbed}
  \begin{split}
    \restrict(M_\epsilon(\tskip;\lift(y)))&
    =\restrict(M_\epsilon(\tskip+\deltat;\lift(x)))\mbox{,}\\
    \restrict(M_\epsilon(\tskip;g_\epsilon(\lift(y_*))))&=
    \restrict(M_\epsilon(\tskip+\deltat;g_\epsilon(\lift(x))))\mbox{,}
  \end{split}
\end{equation}
where the equation at the top defines $y=\Phi(\deltat;x)$. To check
how the difference $y-y_*$ depends on $x$, $\tskip$, $\deltat$, and
$\epsilon$ we use a regular perturbation argument by comparing
solutions of the two equations in \eqref{eq:exactperturbed}. We rely
on \eqref{eq:gepsconv}, which guarantees that the perturbations are
small, and the invertibility of the map \eqref{eq:adef}, which
guarantees that the linearization of the left-hand side with respect
to $y$ and its inverse are uniformly bounded for
$\epsilon(\tskip+\deltat)\leq \tup$.
\begin{theorem}[Convergence of approximate flow map]\label{thm:conv}
  Let $K\in(0,K_0)$ be a given constant.  We assume that the
  assumptions on time scale separation
  \textup{(}Assumption~\ref{ass:timescales}\textup{)} and
  transversality
  \textup{(}Assumption~\ref{ass:transversality}\textup{)} hold for
  $\lift$, $M_\epsilon$ and $\restrict$ such that
  \begin{displaymath}
    x\mapsto\restrict(M_\epsilon(t;g_\epsilon(\lift(x))))
\end{displaymath}
is a local diffeomorphism if $|\epsilon t|\leq\tup$ with some $\tup>0$
that is uniform for all $\epsilon$ and all $x\in\dom\lift$.

Then there exist a lower bound $t_0$ for $\tskip$, an upper bound
$\epsilon_0$ for $\epsilon$, and a constant $C>0$ such that
$y=\Phi(\deltat;x)$ and $y_*=\Phi_*(\deltat;x)$ are well defined by
\eqref{eq:exactperturbed}, and the estimate
  \begin{equation}
    \label{eq:thm:accuracy}
    \|\partial_2^j\Phi(\deltat;x)-\partial_2^j\Phi_*(\deltat;x)\|\leq 
    C\exp(-K\tskip)
  \end{equation}
  holds for all orders $j\in\{0,\ldots,k\}$, all $x\in\dom\lift$,
  $\epsilon\in(0,\epsilon_0)$, $\tskip\in(t_0,\tup/\epsilon]$, 
  and $\deltat\in[0,\tup/\epsilon-\tskip]$.
\end{theorem}

(Remember that $k$ is defined above \eqref{eq:msmall}.)
Note that the assumptions of Theorem~\ref{thm:conv} require that
  $\epsilon\tskip$ and $\epsilon(\tskip+\deltat)$ be bounded by
  $\tup$. Hence, the theorem ensures convergence of $\Phi(\deltat;x)$
  to $\Phi_*(\deltat;x)$ only if $\tskip\to\infty$ and $\epsilon\to0$
  simultaneously. Since $\epsilon$ is usually fixed in applications,
  this theorem is not enough to ensure convergence for
  $\tskip\to\infty$ uniform for $\epsilon$.

The proof of Theorem~\ref{thm:conv} splits the error
$\Phi(\deltat;x)-\Phi_*(\deltat;x)$ using the fiber projection
$g_\epsilon$. The projection of the error onto ${\cal C}_\epsilon$
using $g_\epsilon$ is zero by construction, and the error transversal
to the manifold decays exponentially due to \eqref{eq:gepsconv},
giving a perturbation of order $\exp(-K\tskip)$. This implies that we
can apply the implicit function theorem if
$\epsilon(\tskip+\delta)\leq\tup$, giving an error of order
$\exp(-K\tskip)$ for the difference $y-y_*$ and for the first- and
higher-order derivatives. The details of the proof are given in
Appendix~\ref{sec:app:conv}.
\section{Discussion of the general convergence statement and its
  assumptions}
\label{sec:comp}
Theorem~\ref{thm:conv} is a local statement with respect to $x$,
claiming convergence only in a region $\dom\lift$ in which the
transversality conditions are uniformly satisfied. One has to restrict the
times $\tskip$ and $\deltat$ such that the slow flow
$M_\epsilon(t;g_\epsilon(x))$ cannot leave the region
$g_\epsilon(\dom\restrict)$ for the times $t=\tskip$ and $t=\tskip+\deltat$. This
is appropriate because in many cases, during continuation or
projective integration the maps $\restrict$ and $\lift$ get adapted
(for example, for the traffic problem investigated in Section \ref{sec:eqfree}, $\lift$ is
varied along the curve of macroscopic equilibria).

\subsection{Comparison to the explicit equation-free approach}
The convergence theorem, Theorem~\ref{thm:conv}, implies that for smaller
$\epsilon$ and a longer healing time $\tskip$ the deviation from the
true flow reduces as long as $\epsilon(\tskip+\deltat)\leq\tup$ and
$\epsilon\tskip\leq\tup$. This is in contrast
to 
the approach proposed by \cite{KevrekidisSamaey2009}, where the
coarse flow map  was defined in an explicit way: 
$\Phi_\mathrm{explicit}(\deltat;x)=\restrict(M_\epsilon(\deltat;\lift(x)))$
or $\Phi_\mathrm{explicit}(\deltat;x)=
\restrict(M_\epsilon(\deltat+\tskip;\lift(x)))$
\cite{Makeev2002,Siettos2003,Cisternas2004}. Following this approach,
one would analyze equilibria of the slow flow and their stability by
studying fixed points of the map
\begin{equation}
  \label{eq:eqfree:orig}
  \Phi_\mathrm{explicit}(t;\cdot): x\mapsto  \restrict(M_\epsilon(t;\lift(x)))
\end{equation}
for $x$, where $0<t\ll1/\epsilon$ is chosen such that it includes a
healing time $\tskip$ ($t>\tskip$). 
(Compare \eqref{eq:eqfree:orig} with definition
\eqref{eq:exactperturbed}: $y=\Phi(\delta;x)$ if
$\restrict(M_\epsilon(\tskip;\lift(y)))=\restrict(M_\epsilon(\tskip+\deltat;\lift(x)))$.)
For $1\ll t\ll1/\epsilon$ the map $\Phi_\mathrm{explicit}(t;\cdot)$ is
a perturbation of order $O(\epsilon t)\ll1$ of the map $x\mapsto
\restrict(g_0(\lift(x)))$.  Any flow map on the slow manifold
must be a perturbation of the identity of order $\epsilon t$ for small
$\epsilon t$. Thus, the explicit map $\Phi_\mathrm{explicit}(t;\cdot)$
can be a valid approximation for the flow on the slow manifold in
any coordinates only if $\restrict\circ g_0\circ\lift$ equals the identity
on $\R^d$. Often this requirement is approximated by
$\restrict\circ\lift=\id$, because $g_0$ is in general unknown
\cite{Kevrekidis2004,springerlink:10.1007/s10827-006-3843-z,Li2003,Reppas2010,AIC:AIC690490727}. Note
that there is no $\epsilon$- or $t$ dependence in the limiting map
$\restrict\circ g_0\circ\lift$, resulting in the much more restrictive
condition $\restrict\circ g_0\circ \lift=\id$ than transversality
Assumption~\ref{ass:transversality} on $\restrict$ and
$\lift$. Moreover, $\restrict\circ g_0\circ\lift=\id$ is only a
consistency condition, making it possible for
$\Phi_\mathrm{explicit}(t;\cdot)$ to resemble the map of a slow
flow. If this consistency condition is violated, then
$\Phi_\mathrm{explicit}(t;\cdot)$ will show dynamics independent of
the properties of the flow on the slow manifold. For example, if the
map $\restrict\circ g_0\circ\lift$ has a stable fixed point, then
$\Phi_\mathrm{explicit}(t;\cdot)$ will also have a stable fixed point
independent of the slow flow $M_\epsilon$ on ${\cal C}_\epsilon$.

One way to ensure that the operator
  $\Phi_\mathrm{explicit}$ approximates the slow flow is to construct
  a lifting operator that maps onto the slow manifold ${\cal
  C}_\epsilon$. This has been achieved up to finite order of
$\epsilon$ through constrained-runs corrections to $\lift$
\cite{Zagaris2009,Zagaris2012}. In our notation the first-order
version of this scheme would correspond to defining the lifting
$\lift:\R^d\ni x\mapsto u\in\R^D$ as the (locally unique) $u$
satisfying $\restrict(u)=x$ and $\d/\d t(\restrict^\pitchfork(u))=0$
(zero-derivative principle), where $\restrict^\pitchfork$ is an
arbitrary operator satisfying
$\R^D=\ker\restrict\oplus\ker\restrict^\pitchfork$.  Zagaris \emph{et
  al.}\ \cite{Zagaris2009,Zagaris2012} developed general $m$th-order
versions of this scheme. Vandekerckhove \emph{et al.}\
\cite{Vandekerckhove2011} compared the constrained-runs schemes from
\cite{Zagaris2009,Zagaris2012} to the results of the implicit
expression \eqref{eq:rmatch} (called \textsc{InitMan} in
\cite{Vandekerckhove2011}) for various examples, finding
\eqref{eq:rmatch} uniformly vastly superior in terms of convergence
and performance. Equation \eqref{eq:rmatch} also requires only the
solution of a $d$-dimensional, not a $D$-dimensional, system (usually
$d\ll D$). Using \eqref{eq:phistardef0} it is not necessary
  to find a microscopic state $u$ on the slow manifold matching a
  particular restriction $x$ ($\restrict u = x$). A usage of
  \textsc{InitMan} prefixed at each single step of an explicit
  equation-free scheme would do so and is an alternative. Recognizing that the slow
  flow is given by an implicit ODE from the beginning reduces the
  computational overhead, because matching the restriction is 
  required only at user-specified points.

\subsection{Testing the transversality conditions and choosing the
  healing time and coarse dimension}
The conditions listed in Assumption~\ref{ass:timescales} and
Assumption~\ref{ass:transversality} contain terms that are unknown in
practice. For example, the fiber projection $g_0$ and the tangent
space ${\cal N}_0$ to the slow manifold are both inaccessible because
in many cases one cannot vary the time scale separation parameter
$\epsilon$. However, observing the minimal singular value of the
linearization $\partial_2P_\epsilon(\tskip;x)=\partial/\partial
x[\restrict(M_\epsilon(\tskip;\lift(x)))]$ with respect to $x$ (a
$d$-dimensional matrix) provides an indicator: in points where the
transversality condition is violated, the linearization becomes
singular.

Similarly, the condition number of the linearization
$\partial_2P_\epsilon(\tskip;x)$,
$\cond\partial_2P_\epsilon(\tskip;x)$, guides the choice of the
optimal healing time $\tskip$. All tasks involve solving nonlinear
equations with a Jacobian $\partial_2P_\epsilon(\tskip;x)$. While
the error due to finite time scale separation becomes smaller,
$\cond\partial_2P_\epsilon(\tskip;x)$ can grow with $\tskip$ such
that other errors may become dominant when they are amplified by
$\cond \partial_2P_\epsilon(\tskip;x)$. In particular, when the
microscopic system is a Monte Carlo simulation, a trajectory
$M_\epsilon(t;u)$ is determined via ensemble runs, and the accuracy of
the evaluation of $M_\epsilon$ is only of the order of $1/\sqrt{S}$,
where $S$ is the ensemble size.

The linearization $\partial_2P_\epsilon(\tskip;x)$ also helps to
reveal whether one has too many coarse variables, that is, whether $d$ is too
large such that the flow $M_\epsilon$ restricted to the assumed slow
manifold ${\cal C}_\epsilon$ is not sufficiently slow (still
containing rapidly decaying components). Then
$\partial_2P_\epsilon(\tskip;x)$ becomes close to singular, too. Note
that any solution found, for example, by solving the fixed point
equation \eqref{eq:eqfree:simple} is still a correctly identified
fixed point with correctly identified stability. However, the
linearization of \eqref{eq:eqfree:simple} becomes close to singular.

\subsection{Chaotic and stochastic systems} Barkley, Kevrekidis, and
Stuart \cite{Barkley2006} analyzed how the equation-free approach can be used
to analyze moment maps of stochastic systems or high-dimensional
chaotic systems that converge in a statistical mechanics sense to
low-dimensional stochastic differential equations
(SDEs). These moment maps play then the same role as the
  macroscopic map $\Phi(\deltat;\cdot)$ in our case. The authors of \cite{Barkley2006} observe that
  the choice of $\deltat$ strongly influences the number and stability
  of fixed points. Also the inclusion of additional macroscopic
variables (increasing $d$) changes the results of the equation-free
analysis qualitatively. It is unclear how the implicit scheme
\eqref{eq:phidef} behaves in the situations studied by
\cite{Barkley2006}. While \cite{Barkley2006} also invokes a
separation-of-time-scales argument to study approximation quality for
the stochastic systems, their setting does not fit into the
assumptions underlying Fenichel's theorem but requires weaker notions
of convergence based on averaging over a chaotic attractor (see
\cite{Givon2004} for a review). An adaptation of the analysis in
\cite{Barkley2006}, and possibly further adaptation of the implicit
scheme \eqref{eq:phidef}, is the missing link between
Theorem~\ref{thm:conv} establishing convergence for the idealized
situation, given in Section \ref{sec:thm}, and applications of
equation-free analysis to stochastic or chaotic systems.


\section{Traffic Modeling --- The Optimal Velocity Model}
\label{sec:model}
We now turn to the equation-free analysis of a system that fits into
the framework of implicit equation-free analysis. We will perform some
of the typical tasks listed in Section~\ref{sec:implicit} 
and apply the implicit equation-free analysis introduced in
  Sections~\ref{sec:implicit} and \ref{sec:thm}.

\begin{figure}[t]
\centering
  \includegraphics[width = 0.45 \textwidth]{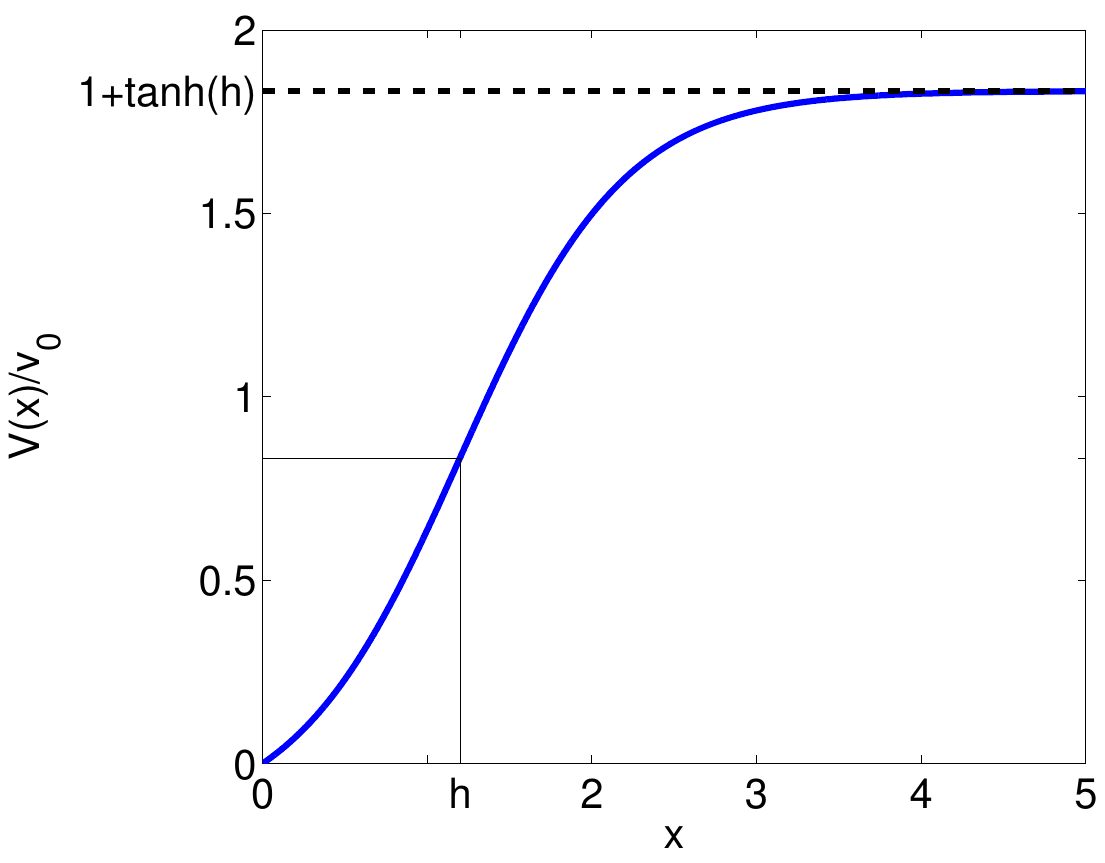}
  \caption{Optimal velocity function $V(x)/v_0$ defined in \eqref{eq:ovfunc}
    for $h=1.2$. The maximal velocity $v_0(1+\tanh(h))$ is obtained
    for $x\to \infty$. $v_0$ acts as a scaling parameter for $V$.
The inflection point of the optimal velocity function $V$ is at $h$.}
  \label{fig:ovfunc}
\end{figure}

We consider $N$ cars driving around a ring road of length $L$. The
individual drivers' behavior is assumed to be uniform and
deterministic, modeled by an optimal velocity model \cite{Bando1995}
of the form
\begin{equation}
  \label{eq:2ode}
  \tau \ddot{x}_n + \dot{x}_n = V(x_{n+1}-x_n), \qquad n = 1,2,\ldots ,N
\end{equation}
where $x_n$ is the position of car $n$, $\tau$ is the inertia of the
driver and car, and $V$ is an optimal velocity function, prescribing the
preferred speed of the driver depending on the distance to the car in
front (the \emph{headway}).  The ring road implies periodic boundary
conditions in space
\begin{equation}
  \label{eq:perbc}
  x_{n+N} = x_n + L.
\end{equation}
In order to do numerical bifurcation analysis, we rewrite the second-order ODE
(\ref{eq:2ode}) as a system of first-order ODEs:
\begin{equation}
\begin{aligned}
  \label{eq:1odesys}
  \dot{x}_n &= y_n\\
\dot{y}_n &= \tau^{-1} \left[ V(x_{n+1} - x_n) - y_n\right].
\end{aligned}
\end{equation}
Similar to \cite{Bando1995,Gaididei2009} we choose the function
\begin{equation}
  \label{eq:ovfunc}
  V(\Delta x_n) = v_0 (\tanh(\Delta x_n-h) + \tanh(h))\mbox{,}
\end{equation}
shown in Figure~\ref{fig:ovfunc}, as the optimal velocity function. In
\eqref{eq:ovfunc}, $v_0(1+\tanh (h))$ is the maximal velocity, $\Delta
x_n := x_{n+1}-x_n$ is the headway, and the inflection point $h$ of
$V$ determines the desired safety distance between cars.
The reviews \cite{Helbing2001,Nagatani2002,Orosz2010} put behavioral
models based on optimal velocity functions into the general context of
traffic modeling and discuss possible choices of optimal velocity
functions. One conclusion from \cite{Orosz2010} is that the choice of
$V$ does not affect the overall bifurcation diagram of a single jam
qualitatively (some choices of $V$ can give rise to unphysical
behavior such as cars briefly moving backwards, though). Depending on
parameters and initial conditions, the system either shows free-flow
behavior, that is, all cars move with the same velocity and headway,
or it develops traffic jams, which means that there coexist regions of
uniformly small headways and low speeds, spatially alternating with regions of
free flow with uniformly large headways and large speeds. 
We focus on the dynamics near the formation of a single jam. In
equilibrium the single traffic jam moves along the ring with nearly
(due to a finite number of cars) constant shape and speed as a
traveling wave against the direction of traffic. In the full system
\eqref{eq:1odesys} the single traffic jam is a traveling wave
perturbed by small periodic oscillations; see
Figure~\ref{sfig:directmacro_b} below.

\subsection{Direct Simulations}
\label{sec:direc}

\begin{figure}[t]
  \centering
  \subfigure[{$h=1.2$,
    $v_0=0.87$}]{\label{sfig:direct_a}\includegraphics[width = 0.45 \textwidth]{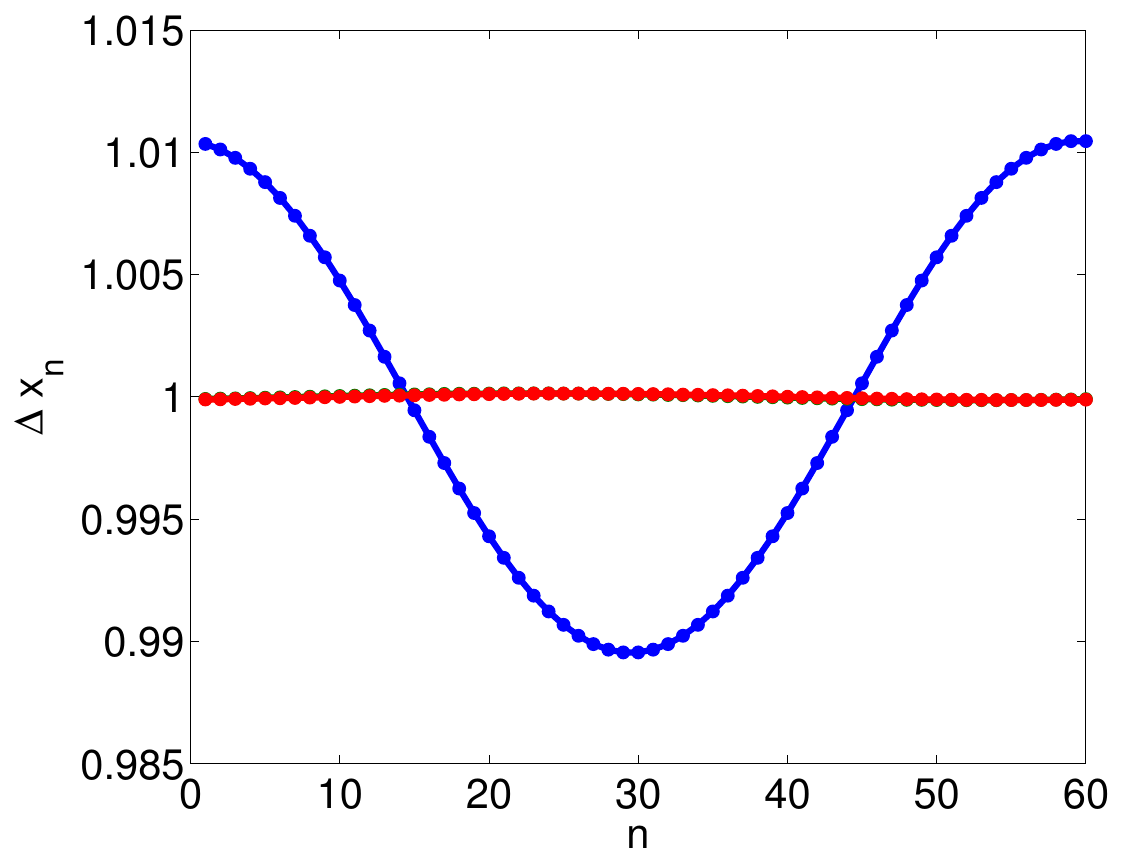}}
\hfill
  \subfigure[{$h=1.2$,
    $v_0=0.91$}]{\label{sfig:direct_b}\includegraphics[width = 0.45 \textwidth]{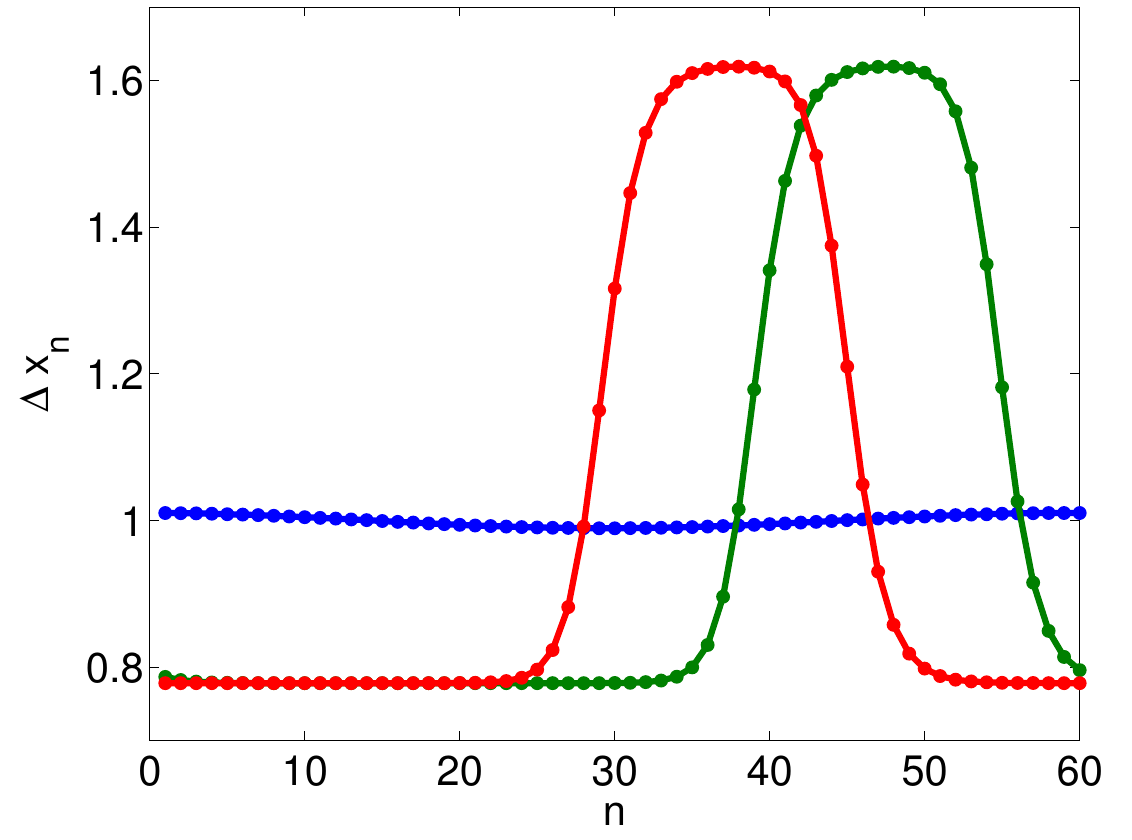}}
  \caption{Comparison of the two traffic flow regimes. The initial
    condition (blue) is compared with the final state \textup{(}red
    $=T=5\times10^4$, green$=T-500$\textup{)}. \textup{(}a\textup{)}
    Free flow regime. \textup{(}b\textup{)} Traffic jam regime. Note the different
    scales of $\Delta x_n$ on the vertical axis.}
  \label{fig:direct}
\end{figure}

\begin{figure}[t]
  \centering \subfigure[{$h=1.2$,
    $v_0=0.87$}]{\label{sfig:directmacro_a}\includegraphics[width =
    0.45 \textwidth]{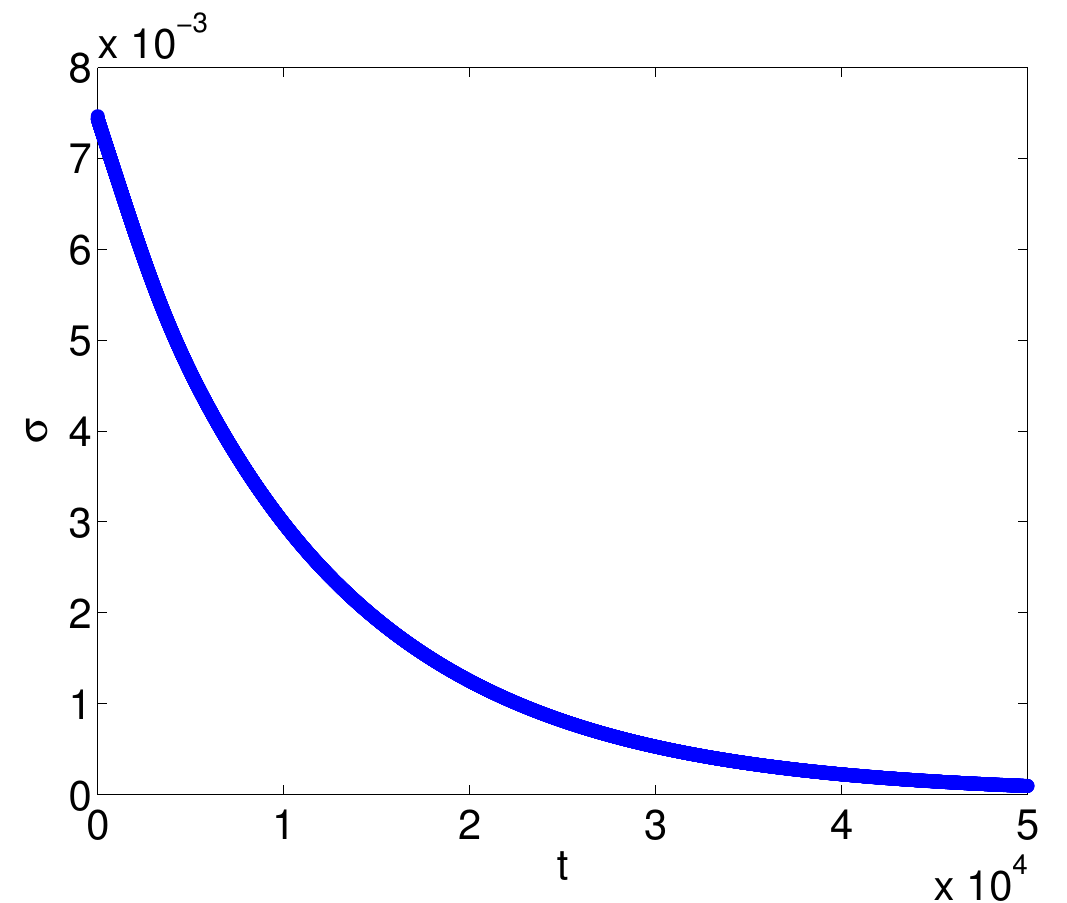}} \hfill
  \subfigure[{$h=1.2$,
    $v_0=0.91$}]{\label{sfig:directmacro_b}\includegraphics[width =
    0.45 \textwidth]{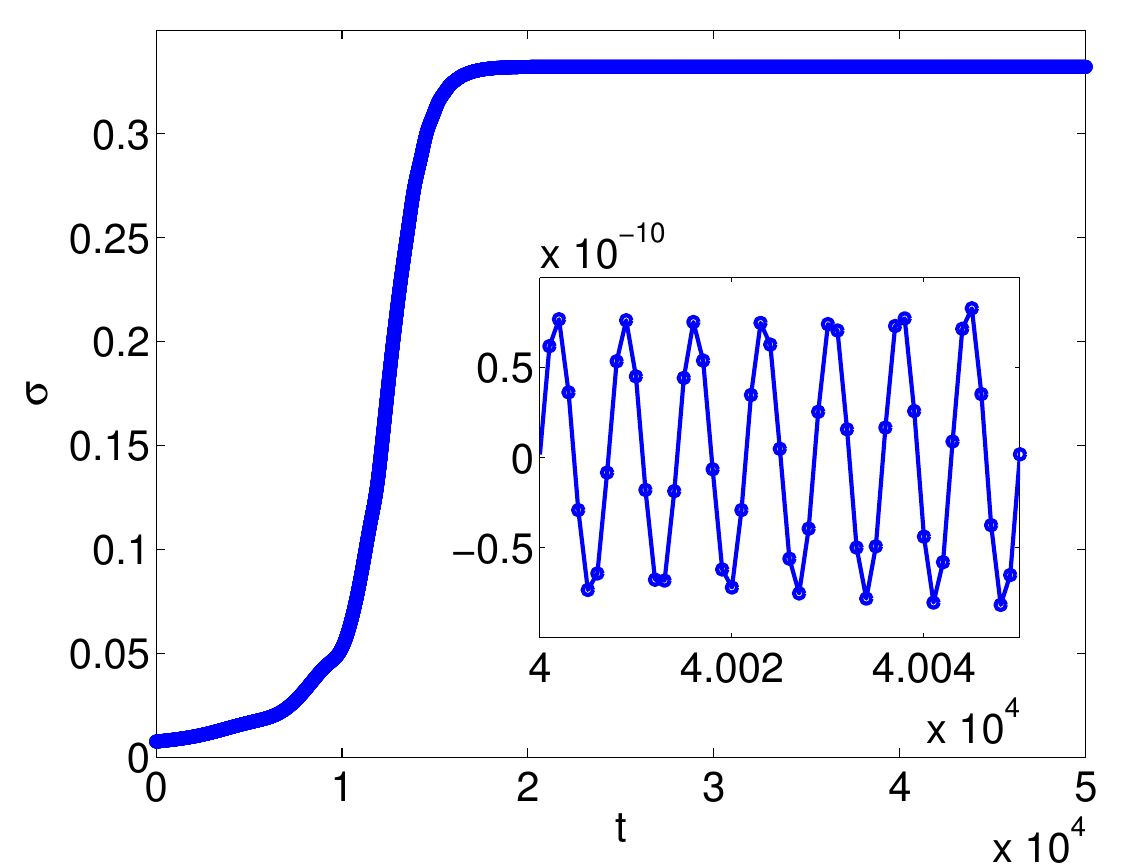}}
  \caption{Time evolution of the macroscopic variable $\sigma$ for the
    same parameters as in Figure \ref{fig:direct}.  \textup{(}a\textup{)} The
    decay to the stable free flow. \textup{(}b\textup{)} Using the same initial
    condition as in \textup{(}a\textup{)}, the system converges to a stable traffic jam.
    The inset in \textup{(}b\textup{)} shows the difference between
    the macroscopic variable $\sigma$ and its long-term average
    $\sigma^*$ over the last $10000$ time steps of the simulation. One expects small oscillations of $\sigma$ in time due
    to the finite number of cars. However, these small oscillations
    are below the tolerance of the ODE
    solver.
  }
  \label{fig:directmacro}
\end{figure}
The uniform flow, starting from initial condition
\begin{equation}
  \begin{aligned}
      \label{eq:freeflow}
      x_n(0) &= (n-1)\frac{L}{N} \\
      y_n(0) &= V\left(\frac{L}{N}\right),
  \end{aligned}
\end{equation}
is a solution of \eqref{eq:1odesys}, where all cars move with the same
velocity $y_n(t)=V(\frac{L}{N})$ and headway $\Delta x_n(t) =
\frac{L}{N}$.  We focus on two types of long-time behavior, the
uniform flow and traveling wave solutions. To give a qualitative
picture of these, we run two simulations, initializing system
\eqref{eq:1odesys} with initial conditions close to the uniform flow,
or adding a periodic perturbation of strength $\mu$:
\begin{equation}
  \begin{aligned}
      \label{eq:ic}
      x_n(0) &= (n-1)\frac{L}{N} + \mu \sin\left(\frac{2\pi}{N} n\right)\\
      y_n(0) &= V\left(\frac{L}{N}\right)\mbox{.}
  \end{aligned}
\end{equation}
For all simulations, we use $N = L =60$.  The simulations were run for
a time $T=5\cdot 10^4$ using the \textsc{Matlab} \texttt{ode45}-solver
\cite{matlab} with absolute and relative tolerance $10^{-8}$. All
parameters for the simulation can be found in table \ref{tab:params}
in Appendix \ref{sec:params}.  For our one-parameter analysis, we also
fix the desired safety distance $h=1.2$. 
Figures~\ref{fig:direct} and \ref{fig:directmacro} show the long-time
behavior of the initial condition \eqref{eq:ic} for the velocity
parameters $v_0 = 0.87$ and $v_0 = 0.91$, respectively.
In Figure \ref{fig:direct}, the headway is shown as a function of car
number. It can be seen that the initial perturbation decays to the
uniform flow for the trajectory for $v_0=0.87$ but converges to a
traveling wave solution for $v_0=0.91$.
\par
We choose the standard deviation $\sigma$ for the headway as the
macroscopic measure (called $x$ in sections \ref{sec:implicit} and
\ref{sec:thm}) describing the traffic flow
\begin{equation}
  \label{eq:macro}
  \sigma = \sqrt{\frac{1}{N-1} \sum_{n=1}^N \left( \Delta x_n -\langle \Delta x \rangle\right)^2}\mbox{,}
  \quad\mbox{where $\Delta x_n=x_{n+1}-x_n$.}
\end{equation}
Here, $\langle \Delta x \rangle = 1$ is the mean of all headways. The
free flow corresponds to $\sigma = 0$ and the decay of $\sigma$ to the
free flow is shown in Figure~\ref{sfig:directmacro_a}. If $v_0$ is
chosen equal to $0.91$, $\sigma$ increases until it settles
to an equilibrium, where a traveling wave of fixed shape is
observed. It can be seen in the inset of Figure
\ref{sfig:directmacro_b} that the macroscopic variable oscillates
even in its steady state. These small-scale oscillations are expected
due to the finite number of cars, because cars arrive at the rear and
leave from the front of the jam at periodic intervals. However, the
oscillation amplitude is orders of magnitude smaller than the
macroscopic dynamics, such that the oscillations are obscured by
discretization effects of the ODE solver (which shows subtolerance
oscillations even for systems with stable equilibria).

\subsection{Time scale separation}
\label{sec:tsep}

\begin{figure}[t]
  \centering
\subfigure[\label{fig:timescale_separation_a}]{\includegraphics[width = 0.45
  \textwidth]{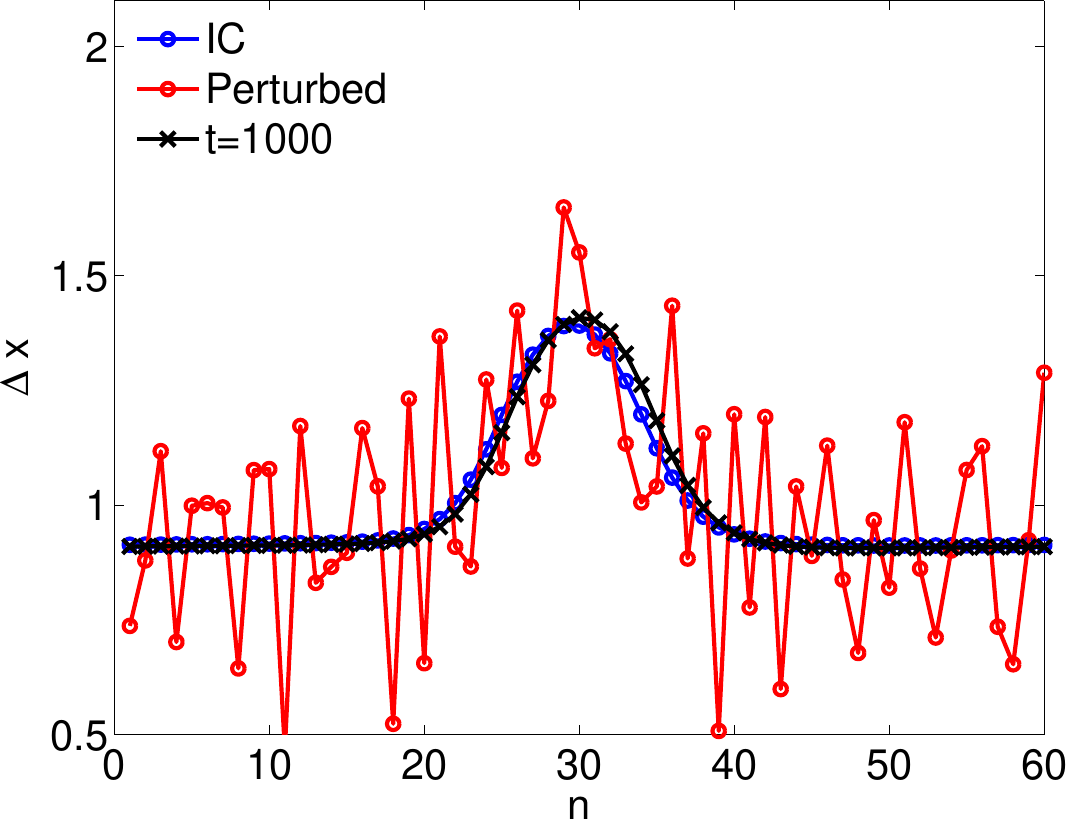}}%
\hfill
\subfigure[\label{fig:timescale_separation_b}]{\includegraphics[width = 0.45 \textwidth]{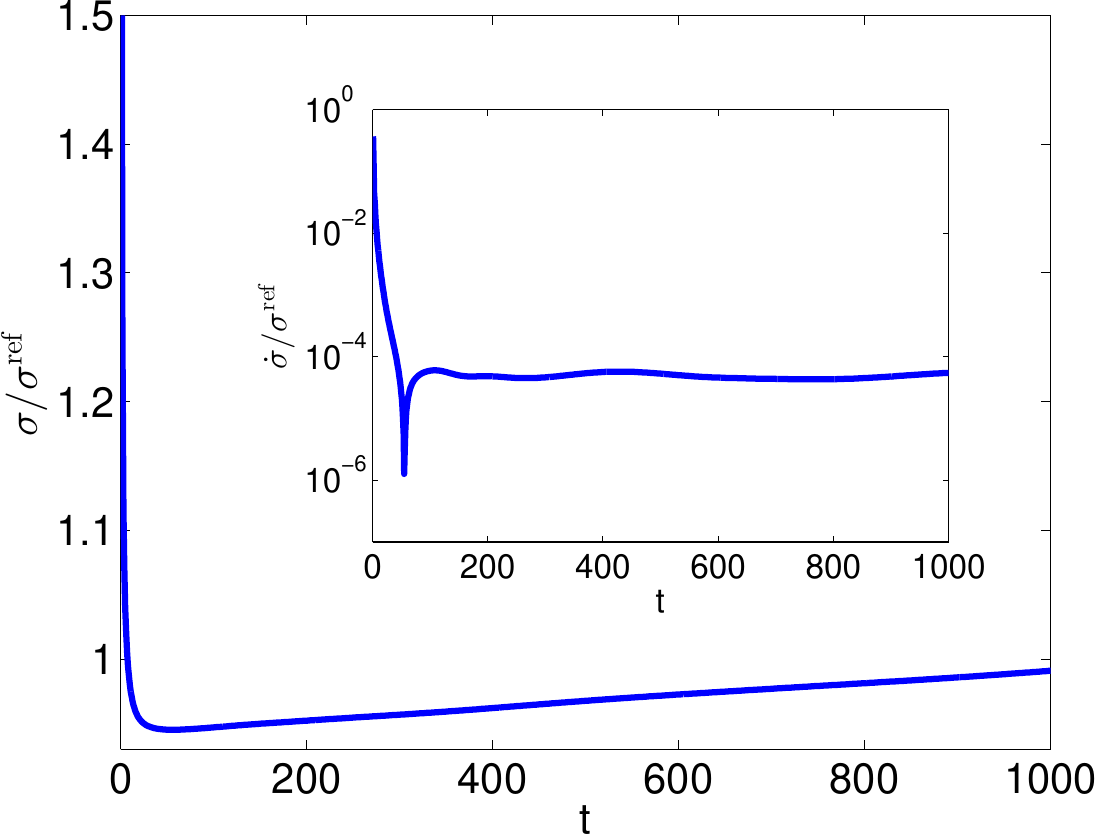}}%
\caption{Visualization of the time scale separation for system
    \eqref{eq:1odesys}. \textup{(}a\textup{)} An initial
    nonequilibrium traffic jam (blue circles) is perturbed with
    uniformly distributed noise to yield a new profile (red
    circles). A microscopic simulation of $1000$ time steps brings the
    system back to a single traffic jam (black crosses), which is
    slightly more pronounced than the initial
    jam. \textup{(}b\textup{)} The same simulation as in
    \textup{(}a\textup{)} shown in the macroscopic variable $\sigma$
    (scaled with $\sigma^\text{ref}=0.16$). After a short time
    ($t\approx 10$), the system relaxes to a one-jam solution. During
    this process, $\sigma$ is decreased drastically and settles on the
    fast time scale (see inset). Afterwards, $\sigma$ increases again
    on a four-orders-of-magnitude slower time scale.} 
  \label{fig:timescale_separation}
\end{figure}
In order to apply the theoretical results from Sections
\ref{sec:implicit} and \ref{sec:thm}, we have to check the extent to which the
assumption about separation of time scales is valid. Initially, we use
simulations to estimate the time scale separation, showing
that the studied one-jam solution forms a one-dimensional stable
submanifold, which we will then study in Section \ref{sec:eqfree}. 

The simulation result shown in Figure \ref{fig:timescale_separation}
highlights that a one-dimensional slow manifold exists 
corresponding to a single jam. 
For Figure \ref{fig:timescale_separation} we perturbed an
initial nonequilibrium traffic jam (blue circles) by adding random
numbers drawn from a uniform distribution in $[-0.5, 0.5]$. This
perturbed state (red circles) is then simulated using
\eqref{eq:1odesys} for $1000$ time steps. The resulting state is
observed to rapidly converge back to a single-jam solution (black
crosses). 
Note that the drift of the jams has been subtracted
in order to center the profiles for a better comparison. In Figure
\ref{fig:timescale_separation_a} the traffic jam at time $t=1000$ is
very slightly more pronounced than the initial jam (which was in
nonequilibrium position, though). The time scale separation can be
observed clearly in the time evolution of the macroscopic variable
$\sigma$ (cf.~Figure \ref{fig:timescale_separation_b}). For a very
short time ($t\approx 10$), the macroscopic variable adjusts
rapidly. This corresponds to the fast scale (see inset in Figure
\ref{fig:timescale_separation_b}). Observing the system for a much
longer time of $t=1000$, the slow drift in the macroscopic variable
corresponds to the slow time scale.  A numerical inspection yields a
time-scale separation of approximately four orders of magnitude, \ie
$\epsilon \approx 10^{-4}$, which appears to be different from $1/N$
(cf.~also Figure \ref{fig:directmacro} for visualizations of the slow
dynamics).

  The next
section presents an equation-free bifurcation analysis for jam
formation on the macroscopic level.

\section{Equation-Free Bifurcation Analysis}
\label{sec:eqfree}
We choose a one-dimensional macroscopic description; that is, the
standard deviation $\sigma$ is the only macroscopic variable.
The change of the chosen macroscopic variable $\sigma$ is studied with
respect to system parameters.  According to the equation-free approach
presented in Section \ref{sec:implicit} the macroscopic ODE has the
implicit form
\begin{equation}
  \label{eq:macrodynsys}
  \frac{\mathrm{d}}{\mathrm{dt}}\restrict(M(\tskip,\lift(\sigma))) = 
  \frac{\partial}{\partial\deltat}
  \restrict(M(\tskip+\deltat,\lift(\sigma)))\Big\vert_{\deltat=0}\mbox{,}
\end{equation}
where the derivative on the right-hand side is approximated by the
finite-difference quotient with finite $\delta$
\begin{equation}
  \label{eq:macrorhs}
  F(\sigma) = \frac{\restrict{M(\tskip+\deltat,
      \lift(\sigma))}-\restrict{M(\tskip,
      \lift(\sigma))}}{\deltat}\mbox{,}
\end{equation}
and $\tskip$ is the healing time, which should be chosen long
  enough for transients to decay (cf.~the discussion in Section \ref{sec:tsep}).

As explained in Section \ref{sec:implicit} and \ref{sec:thm}, the
equation-free setup avoids an analytical derivation of a macroscopic
ODE but uses \eqref{eq:macrodynsys} where \eqref{eq:macrorhs} is
evaluated by simulation bursts of length $\tskip+\deltat$. A good
  choice for the time $\delta$ depends on the slow dynamics. We used
  numerical observations to obtain a good estimate for
  \eqref{eq:macrorhs}, see also Figure
  \ref{fig:timescale_separation}. Note that the left- and right-hand
  sides in \eqref{eq:macrodynsys} depend
also on the system parameters $h$ and $v_0$, which are not expressly
included in \eqref{eq:macrodynsys} and \eqref{eq:macrorhs}. We also
drop the subscript $\epsilon$ of $M$ because it enters our
system only indirectly. 
In order to
find trajectories or equilibria of
\eqref{eq:macrodynsys}--\eqref{eq:macrorhs}, it is necessary to define
a lifting operator $\lift$ and a restriction operator $\restrict$.  In
our case, the restriction operator $\restrict$ is given by the
definition of the macroscopic measure in \eqref{eq:macro}, \ie
\begin{equation}
  \label{eq:restriction}
  \restrict (u) = \sqrt{\frac{1}{N-1} \sum_{n=1}^N \left( \Delta x_n -\langle \Delta x \rangle\right)^2}.
\end{equation}
Our lifting operator constructs initial conditions with the help of a
reference state $\tilde u = (\tilde x, \tilde y)\in\R^{2N}$, obtained
during a previous microscopic simulation.  We have to guaranteet hat
the lifting $\lift$ initializes the system into the vicinity of the
solution of interest, which we described in Section
\ref{sec:thm} as $\lift$ having to map into the attracting
neighborhood ${\cal U}$ of the slow manifold. 

The following description assumes that microscopic simulations start
and end near a single-pulse traffic jam.  The components of the
reference state $\tilde u$ are the positions $(\tilde x_n)_{n=1}^N$
and the velocities $(\tilde y_n)_{n=1}^N$ of the
cars (cf. \eqref{eq:1odesys}).  Let us denote the macroscopic state
corresponding to $\tilde u$ by $\tilde \sigma = \restrict (\tilde
u)$. Given a real parameter $p$, whose meaning we shall explain in
detail below, and a reference state $\tilde u$, we define
$\lift_{p,\tilde u} (\sigma)$ to be
\begin{equation}
  \label{eq:reslift}
\begin{aligned}
\lift_{p,\tilde u} (\sigma) = u=(x,y) &= \left(x_\text{new}, 
  y_\mathrm{new}\right)\in\R^N\times\R^N\mbox{,}
\qquad \text{where} \\
\Delta x_\text{new} &= \frac{p \sigma}{\tilde
  \sigma}\Big(\Delta \tilde x -\avg{\Delta \tilde x}\Big)+\avg{\Delta
  \tilde x}\mbox{,} \\ 
 x_{\text{new},1} = 0 \mbox{,} \qquad x_{\text{new},n} &= \sum_{i=1}^{n-1} \Delta x_{\text{new},i}
 \qquad n=2,\ldots,N \mbox{,}\\
 y_{\mathrm{new},n}&=V(\Delta x_{\text{new},n})
 \qquad n=1,\ldots,N \mbox{.}
\end{aligned}
\end{equation}
$V$ is the optimal velocity function \eqref{eq:ovfunc}, $\avg{\cdot}$
refers to the average of a quantity, and $\Delta \tilde x$ are the
headways of the reference state ($\Delta \tilde x_n=\tilde
x_{n+1}-\tilde x_n$).  In \eqref{eq:reslift} we compute the positions
$x\in\R^N$ first and then initialize the velocities $y\in\R^N$ by
using the optimal velocity function for these positions.  The
positions are initialized such that $x_1 = 0$, resulting in a unique
mapping from headways to positions.  The definition \eqref{eq:reslift}
of $\lift_{p,\tilde u}$ contains an artificial parameter $p$, which we
keep equal to unity throughout, except for
Figure~\ref{fig:lifting_comp} in Section \ref{sec:complift} and
  the error estimates in Section \ref{sec:comptskip}. A
parameter value of $p\neq1$ introduces a systematic bias into our
lifting such that we can vary $p$ gradually to investigate how our
results depend on our choice of lifting. For $p\neq1$, the lifting
$\lift_{p,\tilde u}$ violates the common assumption of equation-free
computations, where the identity $\restrict \circ \lift = \id$ is
claimed to be necessary
\cite{Reppas2010,springerlink:10.1007/s10827-006-3843-z,Kevrekidis2004,Li2003,AIC:AIC690490727}.
An application of $\lift_{p,\tilde u}$ and $\restrict$ without any
time evolution in between, yields $ \restrict (\lift_{p,\tilde u}
(\sigma)) = p \cdot \sigma$.

\par
In the following, we use an equation-free pseudoarclength
continuation scheme to compute bifurcation diagrams for the fixed
point of \eqref{eq:macrodynsys}--\eqref{eq:macrorhs}; that is, we
track a root curve (branch) of
\begin{equation}
  \label{eq:macrofp}
  F(\sigma,v_0) = 0
\end{equation}
in the $(\sigma,v_0)$-plane for the macroscopic right-hand side \eqref{eq:macrorhs}. The influence of speed limits on traffic
jam formation motivates the choice of the velocity parameter $v_0$ as a bifurcation
parameter. In \eqref{eq:macrofp} we include the
bifurcation parameter $v_0$ explicitly as an argument of $F$.  The
pseudoarclength continuation contains two steps. The first step is a
predictor step, where we use a secant predictor, assuming that we know two points on the branch
already.  Let $(\sigma^0, v_0^0)$ and $(\sigma^1, v_0^1)$ be those two
points. We define the secant direction by
\begin{equation}
  \label{eq:tangent}
  w = (\sigma^1 - \sigma^0, v_0^1 - v_0^0).
\end{equation}
The prediction $(\hat \sigma, \hat{v}_0)$ for the next point on the
branch is then determined by the secant predictor 
\begin{equation}
  \label{eq:predictor}
  (\hat \sigma, \hat{v}_0) = (\sigma^1, v_0^1) + s \frac{w}{\norm{w}},
\end{equation}
where we keep the stepsize of the predictor uniformly at
$s=10^{-3}$. The prediction is not exactly on the branch and must be
corrected in the following corrector step, which is chosen to be
perpendicular to the predictor direction \eqref{eq:tangent}. The corrector step
solves the system
\begin{equation}
  \label{eq:corrector}
  \begin{split}
    F(\sigma,v_0) &= 0 \\
    w^{(\sigma)} (\sigma - \hat\sigma) + w^{(v_0)} (v_0 - \hat{v}_0)
    &= 0,
  \end{split}
\end{equation}
where $w^{(\sigma)}$ and $w^{(v_0)}$ are the components of $w$ in the
$\sigma$ and $v_0$ direction, respectively.
System~\eqref{eq:corrector} can be solved with respect to $\sigma$ and
$v_0$ by Newton's method using
\begin{equation}
  \label{eq:newton}
  (\sigma^{k+1},v_0^{k+1})^T = (\sigma^k, v_0^k) + \nu J^{-1} F(\sigma^k,v_0^k),
\end{equation}
where $J$ is the Jacobian of the left-hand side of
\eqref{eq:corrector}, given by
\begin{equation}
  \label{eq:jacobian}
  J =
  \begin{pmatrix}
    F_\sigma & F_{v_0} \\
    w^{(\sigma)} & w^{(v_0)}
  \end{pmatrix}\mbox{,}
\end{equation}
and $\nu$ is a relaxation parameter adjusting the length of a Newton
step.  For all computations we used a full Newton step, that is, $\nu
= 1$.  If the information on the Jacobian of the system is poor, for
example, in noisy or stochastic systems, it might be useful to use a
damped Newton method ($\nu < 1$). The iteration is initialized with
the predictor \eqref{eq:predictor}
\begin{equation}
  \label{eq:newtonini}
  (\sigma^{0},v_0^{0}) = (\hat \sigma, \hat{v}_0)\mbox{.}
\end{equation}
During the iteration the function $F$ has to be evaluated according to
its definition \eqref{eq:macrorhs}. This means that we lift, run the
simulation of the microscopic system and then restrict with
$\tskip=300$ and $\delta=2000$.

The Jacobian $J$ is approximated via finite differences. Since
$w^{(\sigma)}$ and $w^{(v_0)}$ are known from the predictor step, we
only have to determine $F_\sigma$ and $F_{v_0}$.  We evaluated $F$ at
the points
\begin{equation}
  \label{eq:evpoints_1st}
  (\sigma, v_0), \quad   (\sigma+\Delta \sigma, v_0), \quad
  (\sigma, v_0+\Delta v_0)
\end{equation}
and computed the one-sided derivatives
\begin{equation}
  \label{eq:compjac_1st}
  F_\sigma = \frac{F(\sigma + \Delta \sigma,v_0) -F(\sigma,v_0)}{\Delta \sigma}, \qquad
F_{v_0} = \frac{F(\sigma,v_0 + \Delta v_0) - F(\sigma,v_0)}{ \Delta v_0}.
\end{equation}

\begin{figure}[t]
  \centering
  \includegraphics[width = 0.95 \textwidth]{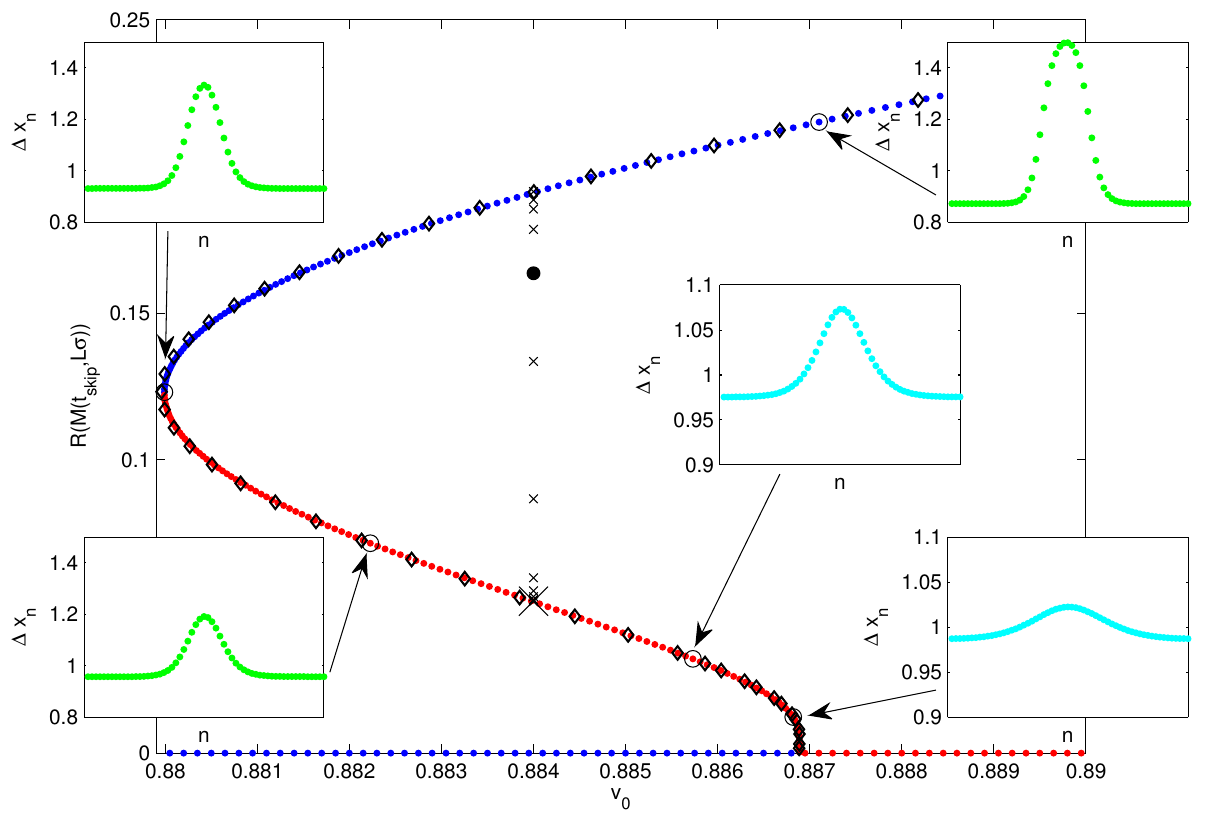}
  \caption{Bifurcation diagram obtained by equation-free
    pseudoarclength continuation for $h=1.2$. The traffic jam
    profiles are shown for selected points in the bifurcation diagram,
    marked with black circles. Note the change in scale on the vertical axes
    on the profiles for better visibility (horizontal axes show the car
    number $n$, vertical axes show headways). A fold point has been detected
    at $(v_0,\restrict(M(\tskip,\lift\sigma))) \approx (0.88,0.125)$,
    where a change in stability is observed. The blue dots mark stable
    states, while the red dots mark unstable states. 
    It is due to the
    equation-free continuation that unstable branches can be
    observed. 
    For the lifting we use \eqref{eq:ic} and \eqref{eq:reslift} for
    continuation of the uniform flow and the traveling wave solution,
    respectively.
    Additionally, the black crosses mark a backward
    trajectory computed by using \eqref{eq:bweuler}. Starting from the
    stable branch, the backward integration converges to the unstable
    branch (big cross). The black dot is the base point used for an
    error estimate in Figure~\ref{fig:tskiperror}. Black diamonds denote
      the results of a direct continuation of the full microscopic system
      on the macroscopic level. The data is in perfect agreement with
      results from implicit equation-free methods.}
  \label{fig:profiles}
\end{figure}

We started the one-parameter continuation of the traffic jam in the
direction of decreasing $v_0$ from two profiles obtained by direct
simulations at $v_0 = 0.91$ and $v_0=0.9$. The resulting bifurcation
diagram is shown in Figure \ref{fig:profiles}.
\par
The traffic jam, \ie traveling wave, is stable for large values of
$v_0$. When following the branch, a saddle-node bifurcation is
detected at $(\restrict(M(\tskip,\lift\sigma))^*,v_0^*) \approx
(0.88,0.125)$, where the traffic jam changes stability. A further
decrease of $v_0$ at that point would make the traffic jam
dissolve. But due to the equation-free pseudoarclength continuation
of the continuous branch, it is possible to follow the branch around
the fold point and continue the unstable branch for increasing
$v_0$. The traffic jam stays unstable until it reaches the uniform
flow at $\sigma = 0$ at a Hopf bifurcation point (cf. Section
\ref{sec:2parcont} and \eqref{eq:v0_analy}).  The microscopic states
corresponding to selected points along the branch are shown as insets
in Figure~\ref{fig:profiles}. The shape has sharp layers and a flat
plateau on the stable branch, and becomes harmonic close to the
equilibrium value $\sigma = 0$.  Additionally, the time steps of a
backward integration are shown for $v_0 = 0.884$, showing the
heteroclinic connection between stable and unstable jams. The
trajectory starts for $t_0=0$ at the stable branch. The Euler
scheme \eqref{eq:expleuler} is used for computing the backward
trajectory; that is,
\begin{equation}
  \label{eq:bweuler}
  \restrict M(\tskip,\lift(\sigma_{j+1})) = 
  \restrict M(\tskip,\lift(\sigma_j)) + F(\sigma_j) \Delta t\mbox{,}
\end{equation}
where $\sigma_j$ is the solution at $t_j=j\Delta t$, and $\Delta t =
-5000$ is chosen. The size of $\Delta t$ is determined by the
  desired accuracy of the coarse projective integration.
For the computation of $F(\sigma)$ the parameters
from Table \ref{tab:params} in Appendix \ref{sec:params} are chosen in \eqref{eq:macrorhs}.
The backward integration converges to the unstable branch. 


\subsection{The influence of the choice of lifting operator}
\label{sec:complift}

\begin{figure}[t]
  \centering
  \subfigure[]{\label{sfig:lifting_comp_a}\includegraphics[width
    = 0.315 \textwidth]{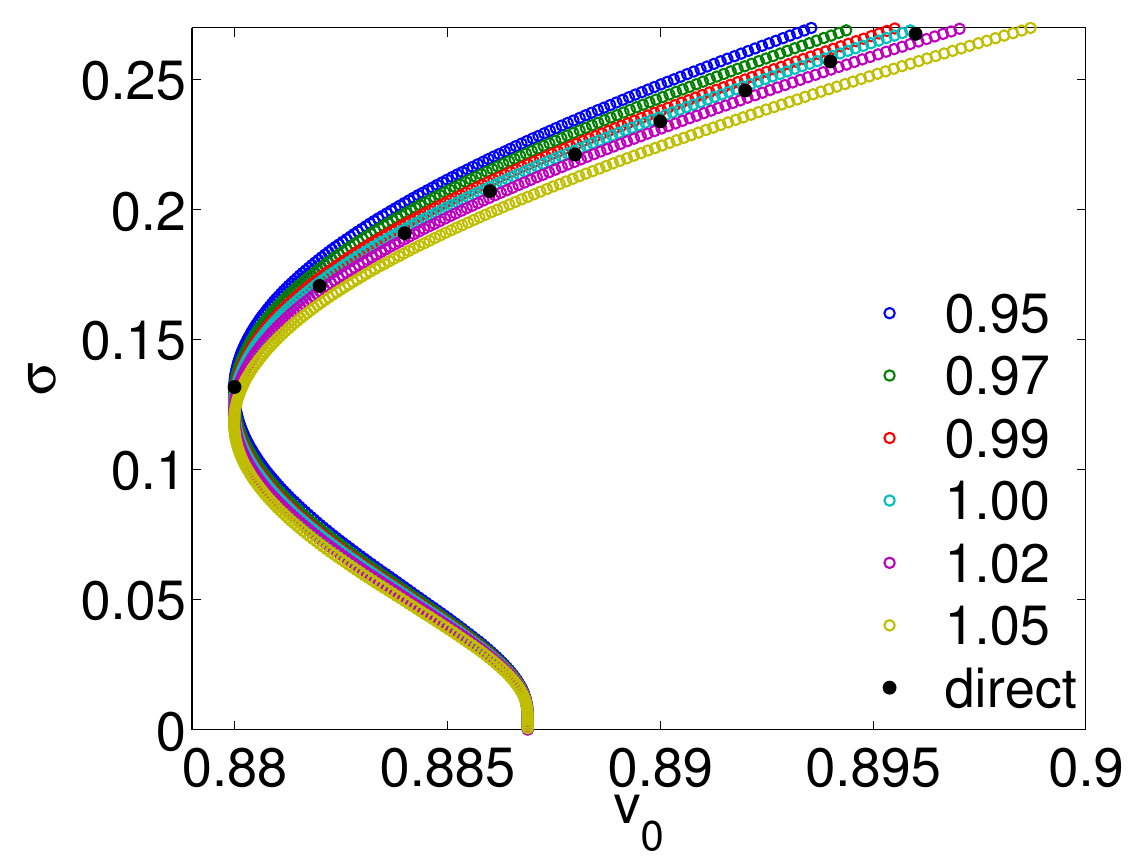}} \hfill
  \subfigure[]{\label{sfig:lifting_comp_b}\includegraphics[width
    = 0.315 \textwidth]{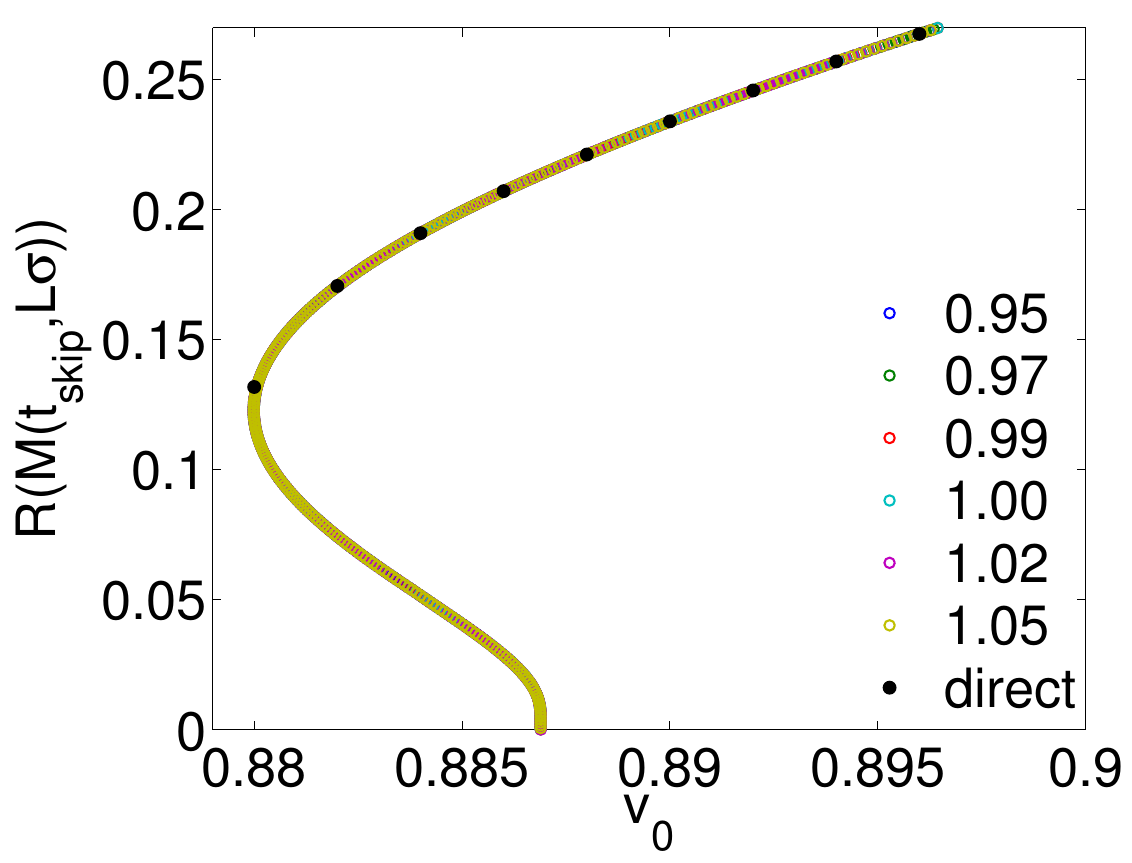}}
  \subfigure[]{\label{sfig:lifting_comp_c}\includegraphics[width
    = 0.315 \textwidth]{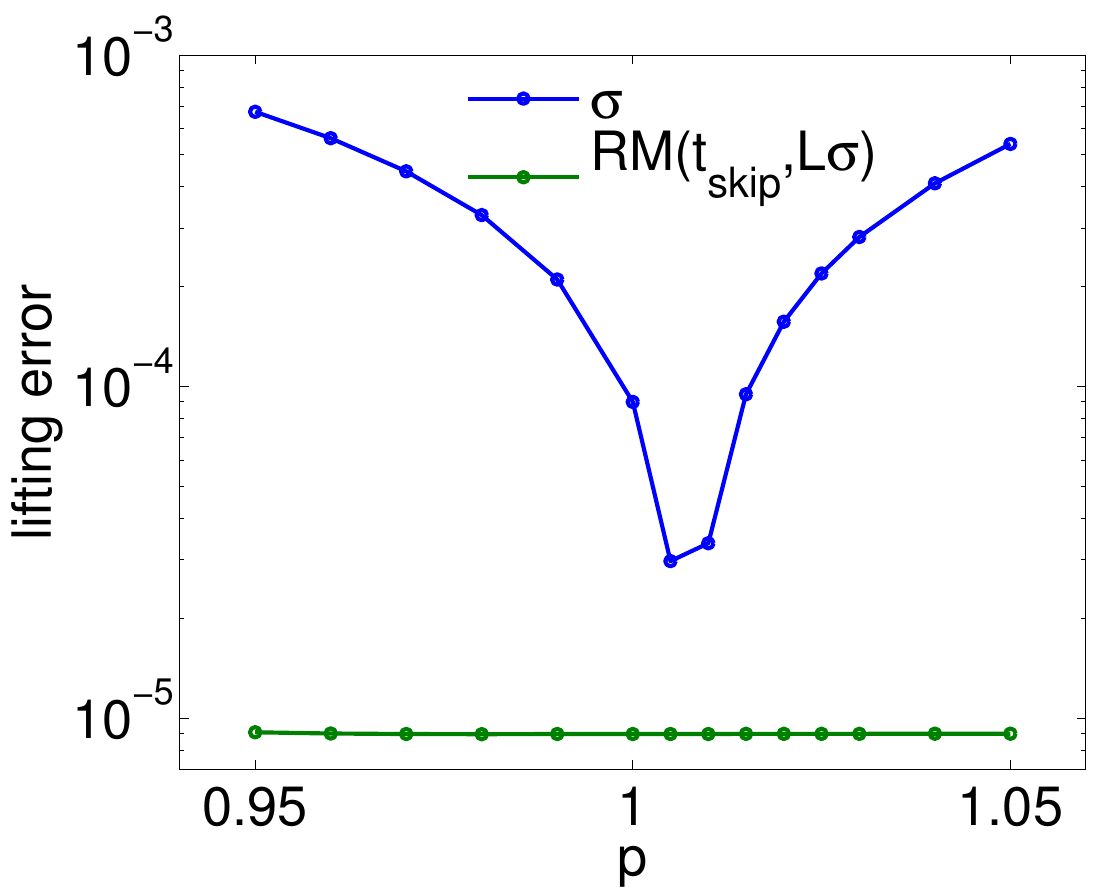}}
  \caption{\textup{(}a\textup{)} Bifurcation diagram obtained from
    running the implicit equation-free continuation scheme described
    in Section \ref{sec:eqfree} for $p=0.95, 0.97, 0.99, 1.0, 1.02,
    1.05$. The black dots show the results from a direct downsweep of
    the stable branch. Depending on the value of $p$, the results
    differ visibly from the direct simulation
    data, which is usually used as an argument for finding a 'good'
    lifting operator. \textup{(}b\textup{)} The healed version of the
    bifurcation diagram; \ie $\sigma_\mathrm{healed} = \restrict
    M(\tskip,\lift_{p,\tilde u} \sigma)$ for different $p$ all
    collapse to the same curve, fitting the direct numerical results
    perfectly. See also the main text in Section
    \ref{sec:implicit}. \textup{(}c\textup{)} Analysis of the lifting
    error. The blue data points show the distances between the
    equation-free solution $v_0(\sigma)$ and the restriction of the
    simulation data using eq. \eqref{eq:l2norm_comp} as a measure for
    the error. The results from a direct simulation of the stable
    branch are used as a reference curve (cf. Figure
    \ref{sfig:lifting_comp_a}). For the 'normal' equation-free data,
    \ie using the unhealed macroscopic quantities, it is observed
    that the error is minimal at $p=1.005$, corresponding to a 'good'
    lifting operator. The green data points show the behavior of the
    healed version of the bifurcation branches (cf. Figure
    \ref{sfig:lifting_comp_b}). The error is uniformly small
    when using the healed data. $h=1.2$ for all images. }
  \label{fig:lifting_comp}
\end{figure}


Figure \ref{fig:lifting_comp} shows how the results depend on the
artificial parameter $p$, which we introduced into the lifting
operator $\lift_{p,\tilde u}$. In both panels, the same bifurcation diagram is
shown for several values of $p$ and compared to the restrictions of
the stable fixed points of direct long-time simulations ($T=3\cdot
10^5$, black dots). The case where the usual equation-free identity
$\restrict \circ \lift_{p,\tilde u} = I$ is fulfilled corresponds to $p=1$. We
observe that the preimages $\sigma$ of the equilibria under the
combination of lifting operator and healing $M(\tskip;\lift_{p,\tilde u}(\cdot))$
depend visibly on $p$ (panel (a) of
Figure~\ref{fig:lifting_comp}). 
Therefore, we compare Figure~\ref{sfig:lifting_comp_a}
with the corresponding Figure~\ref{sfig:lifting_comp_b} for the healed
macroscopic quantity
\begin{equation}
  \label{eq:healedmacro}
  \sigma_\mathrm{healed}=\restrict(M(\tskip,\lift_{p,\tilde u} \sigma))
\end{equation}
for each macroscopic equilibrium $\sigma$ along the branch of the
bifurcation diagram. According
to Section~\ref{sec:thm} the map $\restrict(M(\tskip,\lift_{p,\tilde u} \sigma))$
is a local diffeomorphism from $\R^d$ into $\R^d$ with $d=1$. Plotting the
bifurcation diagram in the $(v_0,\sigma_\mathrm{healed})$-plane in
Figure~\ref{sfig:lifting_comp_b}, we obtain a solution branch that is
independent of the choice of the lifting operator, as
one would expect from Theorem~\ref{thm:conv}.  

For a more detailed analysis of the error, we compute the
$\mathbb{L}_2$ norm between the interpolated data sets $v_0(\sigma)$
(expressing the parameter as a function of the equilibrium location
near the fold) for the direct simulation data and the data for the
stable branch of the equation-free bifurcation diagram. For
interpolation, the \textsc{Matlab} \texttt{interp1} function \cite{matlab} with
the ``\texttt{spline}'' option is used.  We use the error measure
\begin{equation}
  \label{eq:l2norm_comp}
  \norm{f - g}^2 = \int_a^b [f(\sigma) - g(\sigma)]^2 d \sigma
\end{equation}
to analyze the deviation between the restriction of the direct
simulation data and equation-free continuation data. Here, $f$ and $g$
are the interpolated data sets $v_0(\sigma)$ for the simulated data
and the equation-free data, respectively, in the range of $\sigma$
between $a=0.125$ and $b=0.25$. 
The unstable branches cannot be compared with direct integration of
the system.  The deviation $E$ using eq. \eqref{eq:l2norm_comp}) with
lifting parameter $p$ is shown in
Figure~\ref{sfig:lifting_comp_c}. The blue data points correspond to
the distance between the restriction of the simulation data and the
equation-free solutions (that is, the preimages of the equation-free
microscopic solutions under $M(\tskip;\lift_{p,\tilde u}(\cdot))$ in
the domain of $\lift_{p,\tilde u}$) . The distance is small for values
of $p$ close to $1$, where the usual identity $\restrict \circ \lift =
I$ is fulfilled. However, the distance for $\sigma_\mathrm{healed}$
(green data) is uniformly small, independent of the choice of
$p$. Therefore, healed quantities should be used when comparing
equation-free results to restrictions of the direct
simulation data. The uniformly small errors in
Figure~\ref{sfig:lifting_comp_c} (in green) suggest that with implicit
time steppers the results are not sensitive to the choice of the
lifting operator. This is in contrast to most equation-free
applications \cite{KevrekidisSamaey2009,Cisternas2004,Kevrekidis2004},
which use explicit time steppers of the form
$\Phi(\delta;x)=\restrict(M(\delta;\lift(x)))$.

\subsection{Influence of the healing time $\tskip$ and comparison
  to explicit scheme}
\label{sec:comptskip}
\begin{figure}[t]
  \centering \subfigure[\label{fig:tskipbif_a}]{\includegraphics[width = 0.45
    \textwidth]{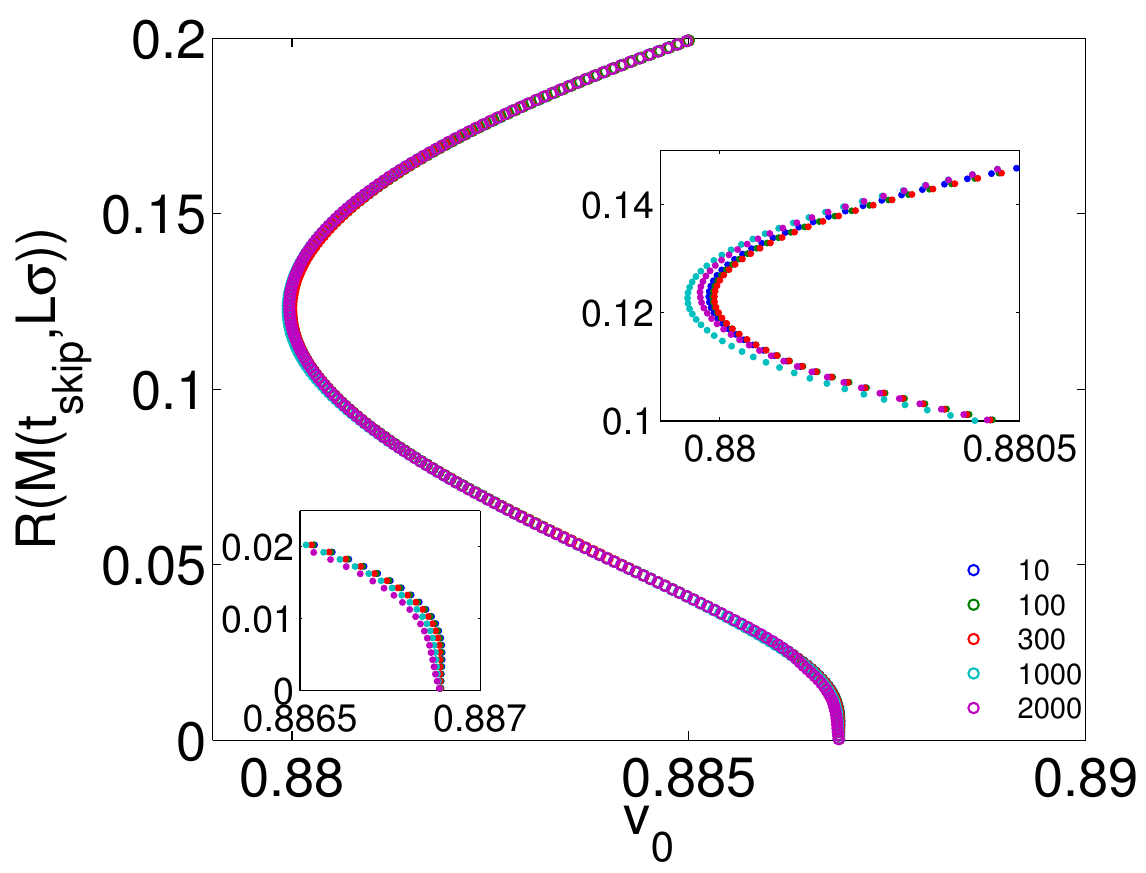}} \hfill
  \subfigure[\label{fig:tskipbif_b}]{\includegraphics[width = 0.45
    \textwidth]{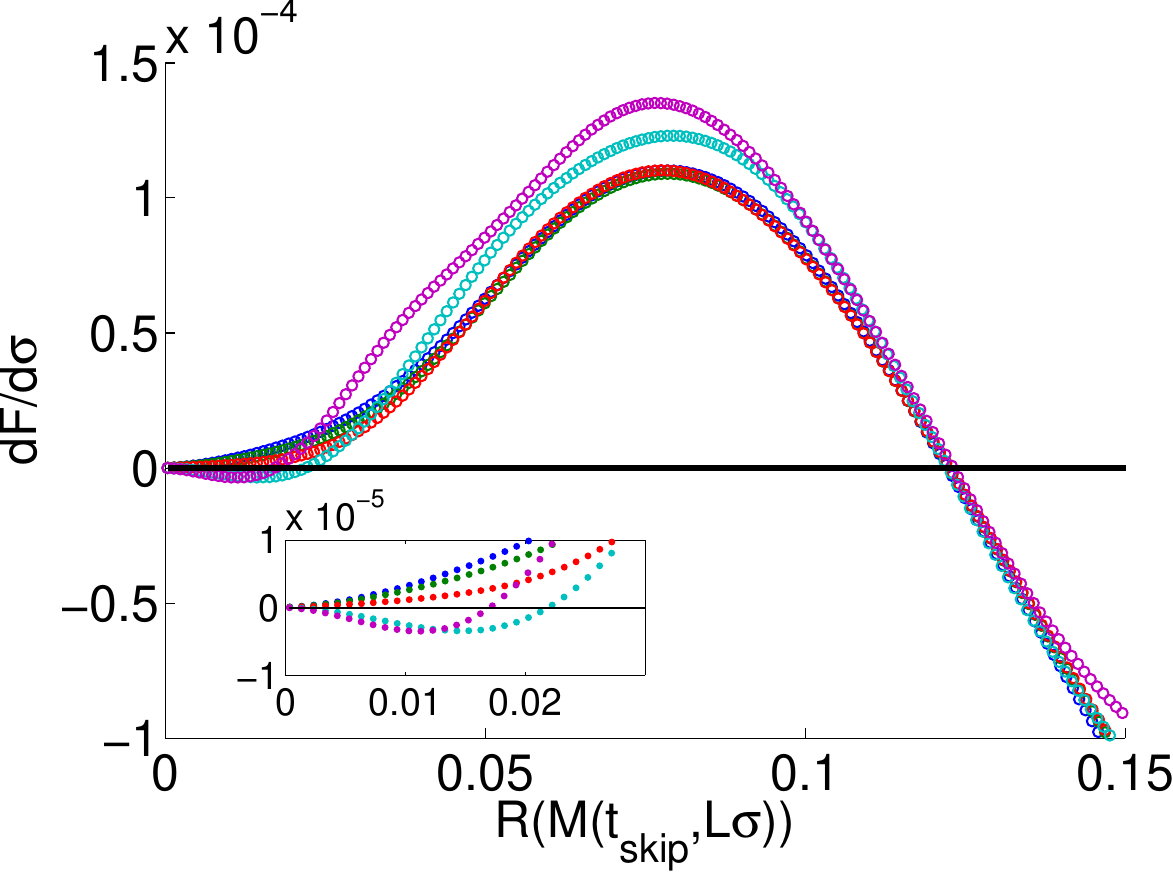}} 
  \caption{\textup{(}a\textup{)} Bifurcation diagrams for $h = 1.2$ and $\tskip =
      10,100,300,1000,2000$ in the healed quantities of $\sigma$. The difference between the curves is
    very small. Insets show a zoom for the fold and the Hopf
    point. \textup{(}b\textup{)} Comparison of the Jacobians for the
    different values of $\tskip$ along the curve. Close to the Hopf point
    the value of the Jacobian does not converge for increasing $\tskip$ within the plotted range.%
  }
  \label{fig:tskipbif}
\end{figure}

In this section, we investigate the influence of $\tskip$ on the
equation-free results, \eg bifurcation diagrams and stability
analysis. First, we show that the bifurcation diagrams are rather
insensitive to the choice of $\tskip$, while the information of the
Jacobian depends more noticeably on the value of $\tskip$.

The bifurcation diagrams obtained for $h=1.2$ and
$\tskip =10,100,300,1000,2000$ are
shown in Figure \ref{fig:tskipbif}. In Figure \ref{fig:tskipbif_a} it
can be observed that the bifurcation diagrams are similar for all
choices of $\tskip$; \ie they show the same qualitative
features. Although the bifurcation diagrams are quantitatively close
to each other, the information about the derivatives, \ie the Jacobian
$\partial F/\partial\sigma$, does not appear to converge to a
  particular value for increasing $\tskip$. Note that the
  scale of the $y$-axis is $10^{-4}$ in this region of the bifurcation
  diagram, which suggests that the slow time scale is of this order.

 For all values of $\tskip$, the fold point near $\sigma=0.12$
is detected by a sign change in the Jacobian (cf. Figure
\ref{fig:tskipbif_b}). Close to the Hopf point (which would appear as
a pitchfork bifurcation in the macroscopic system
\eqref{eq:macrodynsys}) the derivative $\partial F/\partial\sigma$ is
not sufficiently accurate to resolve the criticality of the Hopf
(pitchfork) bifurcation, which appears to be close to being degenerate.  The
Hopf bifurcation point cannot be studied using the operator
$\restrict$ because expression \eqref{eq:restriction}, defining
$\restrict$, is singular in the uniform flow. 

\begin{figure}[t]
  \centering
\includegraphics[width = 0.60 \textwidth]{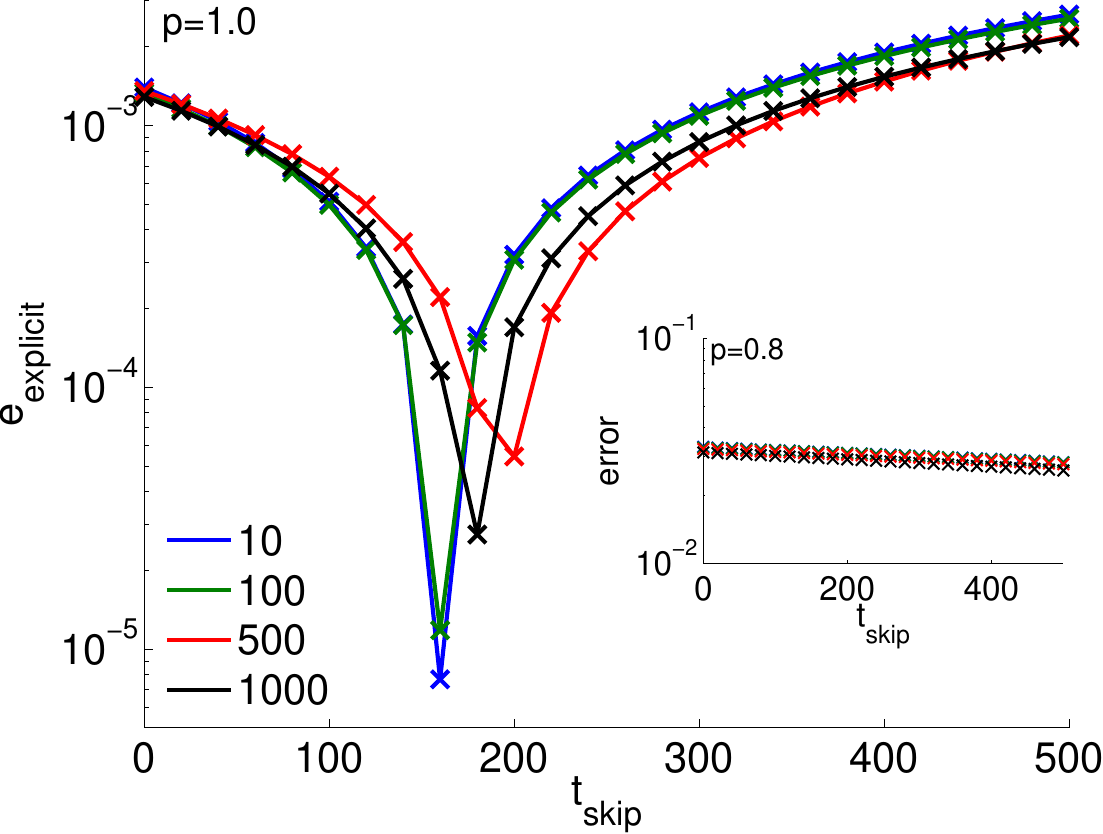}
\caption{Error analysis for the healing time $\tskip$ in the
    explicit scheme, showing the error
    $e_\mathrm{explicit}=|\Phi_\mathrm{explicit}-\Phi_*|$ defined in
    \eqref{eq:explerror}. Colors indicate different values for
    $\delta$. The inset shows the same computation for a scaled
    lifting operator $\lift_{p,\tilde u}$ with $p=0.8$. Here, the
    explicit method has an error which is about two orders of
    magnitude larger than that for a good lifting operator
    $p=1.0$. Note that explicit equation-free computations usually
    require $\restrict \circ \lift = \id$, and the choice of $p=0.8$
    violates this assumption.
  }
  \label{fig:tskiperror}
\end{figure}

\begin{figure}[t]
  \centering
  \subfigure[]{\label{sfig:tskipmeth2}\includegraphics[width = 0.45
    \textwidth]{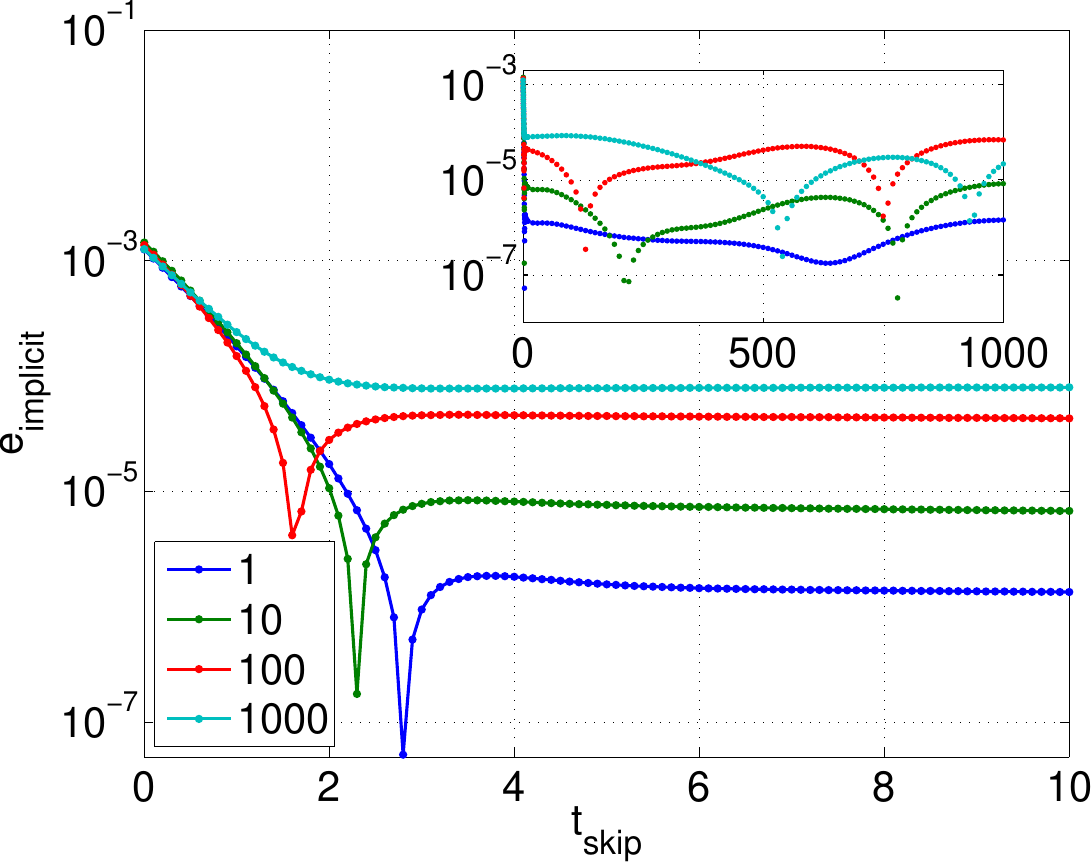}} \hfill
  \subfigure[]{\label{sfig:floquetosc}\includegraphics[width = 0.45
    \textwidth]{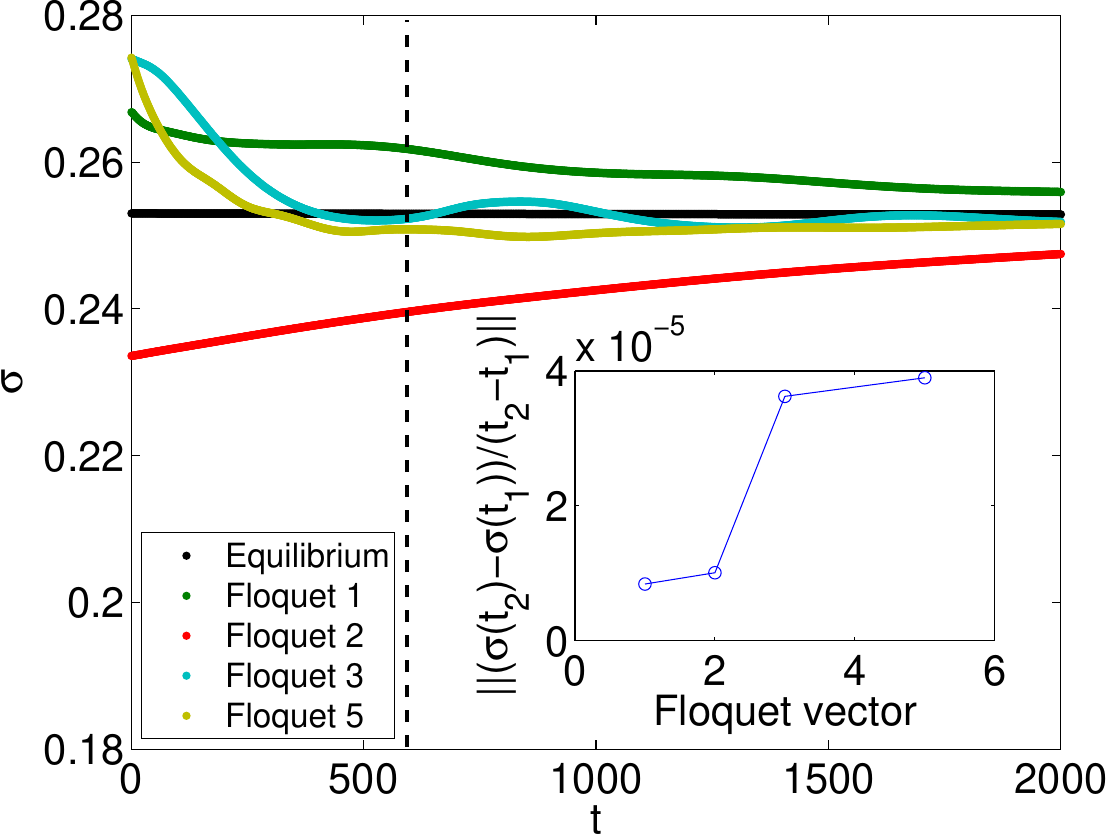}}
  \caption{\textup{(}a\textup{)} Dependence of the error
      $e_\mathrm{implicit}$ given by \eqref{eq:tskiperrmeth2}.  The
      error is shown in dependence on $\tskip$ for $\delta =
      1,10,100,1000$ (see color code in the legend); $h=1.2, v_0=0.884, p=0.8$. 
      \textup{(}b\textup{)} Evolution of the system for perturbations
      in the directions of the leading Floquet vectors (cf. also
      Figure \ref{fig:floquet} for leading Floquet exponents). These
      perturbations lead to oscillations in the macroscopic
      description. 
      The inset shows the decay rate over a time $t_2-t_1 =
      600$. While the perturbations in the first two Floquet
      eigenvectors decay with almost the same rate, there is a gap to
      the next Floquet vector number $3$.}
  \label{fig:errmeth2}
\end{figure}

To study the influence of $\tskip$ on the explicit scheme
$\Phi_\mathrm{explicit}$ and the implicit scheme $\Phi$ in more detail, we compare
the results generated by the approximate macroscopic flow directly to
a pregenerated trajectory of the microscopic flow. To this end we
perform a long-term microscopic simulation of the traffic model from a
reference point (Figure \ref{fig:profiles}, black dot). After a
sufficiently long transient, the dynamics settle to the slow
manifold. We denote the point at the end of this transient as $u(0) =
\tilde u$. The microscopic trajectory $u(t)$ starting from $\tilde u$
is always close to the slow manifold.  The macroscopic state
corresponding to $\tilde u$ is denoted by $\tilde \sigma =
\restrict(\tilde u)$.  The error of the explicit equation-free
approach (scheme \eqref{eq:eqfree:orig}) is then
\begin{equation}
  \label{eq:explerror}
  \begin{split}
    e_\mathrm{explicit}(\tskip,\delta;\restrict(\tilde u))&=
    |\Phi_\mathrm{explicit}(\delta; \restrict(\tilde
      u))-\restrict(u(\delta))|\\
    &= |\restrict(M(\tskip+\delta;\lift_{p,\tilde
        u}\restrict(\tilde u)))-\restrict u(\delta)|\mbox{.}
  \end{split}
\end{equation}
Figure~\ref{fig:tskiperror} shows this error for several fixed
$\delta$ and varying $\tskip$.



The error $e_\mathrm{explicit}$ is of order $10^{-3}$ to $10^{-5}$ for
a good lifting operator, \ie $p=1.0$.  The downward peak around
$\tskip \approx 150$ in Figure~\ref{fig:tskiperror} in logarithmic
scale corresponds to a sign change of the scalar quantity
$\restrict(M(\tskip+\delta;\lift_{p,\tilde u}\restrict(\tilde
u)))-\restrict u(\delta)$ in \eqref{eq:explerror}. For this healing
time $\tskip\approx 150$ the lifted state is mapped into the stable
fiber corresponding to $\tilde u$; that is, $\tilde u=
g_\epsilon(M(\tskip;\lift_{p,\tilde u}(\restrict(\tilde u))))$. Note
that for a one-dimensional slow manifold the stable fibers are
codimension-one surfaces (called isochrones if the slow manifold is a
periodic orbit) such that we can expect to find the fiber for which
the error goes to zero for $\delta\to\infty$ by varying the healing
time $\tskip$.  However, this appropriate healing time may depend on
the point $\tilde u$ on the slow manifold and is in general not
known. The inset in Figure \ref{fig:tskiperror} shows the error
$e_\mathrm{explicit}$ for a nonoptimal lifting operator
$\lift_{p,\tilde u}$, namely for $p=0.8$. The error for the explicit
method is of order $10^{-1}$ to $10^{-2}$ uniformly for $\tskip$ and
$\delta$. Hence, for the explicit scheme varying $\tskip$ can in
general not compensate for errors introduced by the lifting operator.

When estimating the error $e_\mathrm{implicit}$ of the implicit
scheme we have to first find the point $\sigma$ corresponding to
$\tilde \sigma$ after healing. Hence, the error $e_\mathrm{implicit}$ is given as
\begin{equation}
  \label{eq:tskiperrmeth2}
  \begin{aligned}
    e_\mathrm{implicit}(\tskip,\delta;\restrict(\tilde u))=&
    |\restrict(M(\tskip+\delta;\lift_{p,\tilde u}(\sigma)))-\restrict(u(\delta))| 
    &&\mbox{where $\sigma$ solves}\\
    \restrict(\tilde u)=&\restrict(M(\tskip;\lift_{p,\tilde
      u}(\sigma)))\mbox{.}
  \end{aligned}
\end{equation}
Figure~\ref{fig:errmeth2}(a) shows $e_\mathrm{implicit}$ for $p=0.8$
(such that the lifting operator is expected to be at some distance from
the slow manifold initially), the same fixed integration times
$\delta$ as in Figure~\ref{fig:tskiperror}, and a range of $\tskip$
from $0$ to $1000$ (see inset in Figure~\ref{fig:errmeth2}(a)).

After an initial decay over a few orders of magnitude (see Figure
\ref{sfig:tskipmeth2} main graph) the error starts to oscillate (see
inset in Figure \ref{sfig:tskipmeth2}) on a small scale compared to
the value of the macroscopic variable.
\begin{figure}[t]
  \centering \subfigure[Floquet exponents]{\label{fig:floquet:mult}\includegraphics[width = 0.45
    \textwidth]{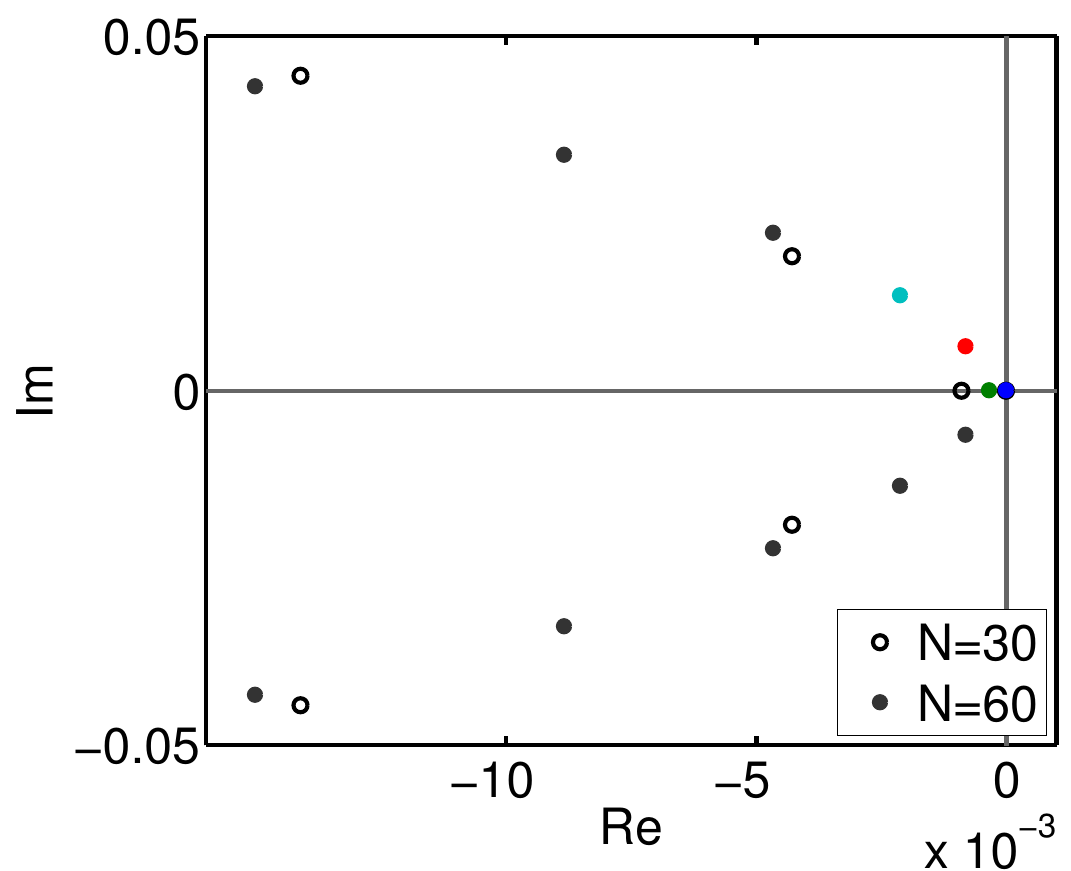}}%
  \hfill \subfigure[Perturbations]{\label{fig:floquet:vec}\includegraphics[width = 0.45
    \textwidth]{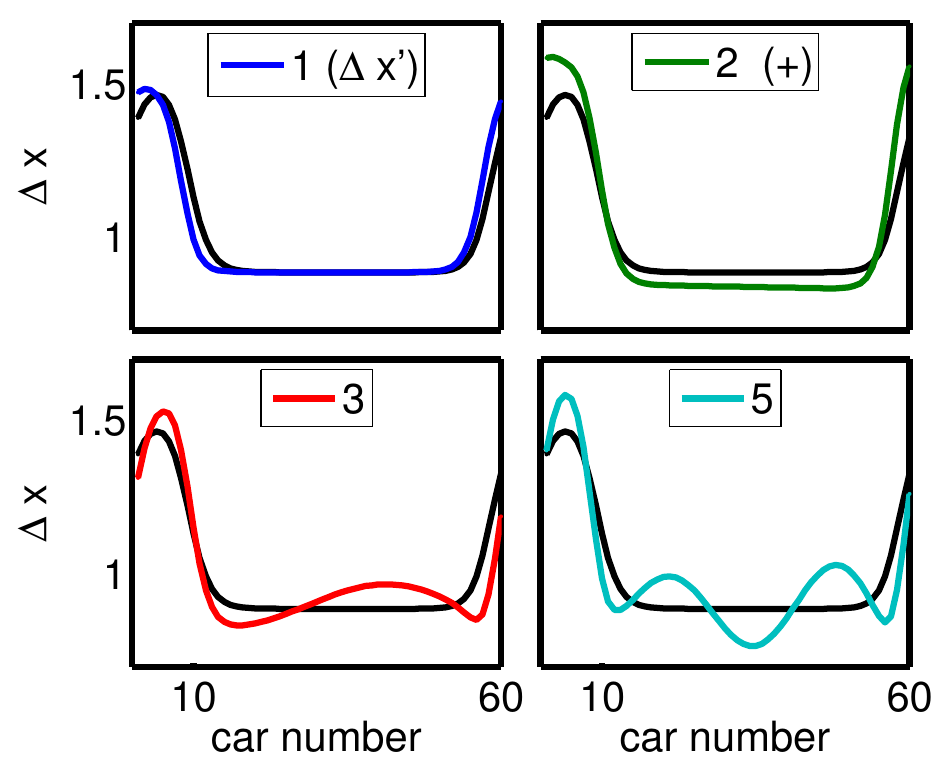}}%
  \caption{{Floquet exponents \textup{(}a\textup{)} and leading Floquet vectors
        \textup{(}b\textup{)} of a single traffic jam \textup{(}viewed
        as a periodic orbit of the full system
        \eqref{eq:1odesys}}\textup{)}. The orbit (\textup{}shown in
      \textup{(}b\textup{))} is also highlighted in
      Figure~\ref{fig:profiles} ($v_0=0.884$ on stable branch). In
      (a), we have included the spectrum of a comparable orbit for $N=30$,
      $L=30$. In \textup{(}b\textup{)} we have added the Floquet vectors
      for the dominant Floquet exponents as a perturbation to
      the periodic orbit. Vectors 3 and 5 are complex. \textup{(+)}:
      Vector 2 has been orthonormalized with respect to vector 1.}
  \label{fig:floquet}
\end{figure}
These small-scale oscillations suggest that the assumptions of
Theorem~\ref{thm:conv} on large time scale separation are not
satisfied for the traffic flow $M$. To confirm this we compute the
Floquet exponents for the stable stationary single-traffic-jam
solution (diamond at $v_0=0.884$ at the end of the heteroclinic
connection marked by crosses in Figure~\ref{fig:profiles}). This is a
periodic orbit of the microscopic system \eqref{eq:1odesys}. 
Figure~\ref{fig:floquet:mult} shows the leading Floquet exponents for
this periodic orbit. It shows a dominant real Floquet exponent very
close to the origin next to the trivial Floquet exponent $0$ (which
corresponds to the flow direction). This dominant real Floquet
exponent corresponds to the slow time scale that the equation-free
analysis attempts to capture. 

Figure~\ref{fig:floquet:mult} also shows
that this dominant Floquet exponent is part of a band of complex
Floquet exponents that is parabola-shaped and bending toward the
half-plane with negative real part (see, for example, the band of full
dots in Figure~\ref{fig:floquet:mult}). The spectra for the two system
sizes plotted in Figure~\ref{fig:floquet:mult} indicate that the
spacing of the Floquet exponents' frequency decreases with increasing $N$. The
parabolic shape of the band then gives a gradually increasing spectral
gap for the low-frequency Floquet exponents until finite-size effects
become visible (to the right of the part of the complex plane shown in
Figure~\ref{fig:floquet:mult}). The spectral gap between the dominant
and the following Floquet exponents gives an upper bound on the time
scale separation that is much more restrictive than the initial
assessment in Figure~\ref{fig:timescale_separation} suggested. 

An
explanation for the apparent discrepancy is the mode shape of the
eigenvectors corresponding to the low-frequency (slow-decay) Floquet
exponents shown in
Figure~\ref{fig:floquet:vec}. Figure~\ref{fig:floquet:vec} illustrates
how perturbations into the directions of the eigenvectors for the
first five Floquet exponents look (ordered by descending real
parts of the exponent). The first Floquet vector corresponds to the
time derivative (the linearization of the time shift). The second
Floquet vector corresponds to the dominant real exponent, tangent to
the slow manifold that the equation-free approach tries to
capture. Floquet vector $2$ is shown orthonormalized with respect to
Floquet vector $1$, because both Floquet vectors $1$ and $2$ are
nearly linearly dependent. While Floquet vector $2$ corresponds to a
change of amplitude of the shape of the jam, the complex Floquet
vectors correspond to spatial perturbations of the jam of low
frequency (the spatial frequency is increasing with increasing time
frequency and decay rate). When decomposing the perturbation given in
Figure~\ref{fig:timescale_separation_a} into the eigenbasis, the
contribution of the space corresponding to the low-frequency,
slow-decay Floquet vectors was small such that one can observe only
small-amplitude low-frequency oscillations after the initial rapid
decay of all high-frequency strong-decay directions (see inset in
Figure~\ref{fig:timescale_separation_b}).

These results explain the oscillations observed in
Figure~\ref{sfig:tskipmeth2}. A perturbation of an equilibrium traffic
jam in the directions of the leading Floquet vectors is shown in
Figure \ref{sfig:floquetosc}. Small-scale oscillations are visible in
the macroscopic trajectories. These oscillations lead to additional
oscillations in Figure \ref{sfig:tskipmeth2} after an initial rapid
exponential decay of the error. Consequently, Theorem~\ref{thm:conv}
is, strictly-speaking, valid only up to a small
residual, 
which in our system is much smaller than the overall dynamics. Thus,
the equation-free approach is applicable (and implicit schemes have
smaller error than explicit ones) even if the conditions of
Theorem~\ref{thm:conv} are not met.

\subsection{Continuation of the fold in two parameters}
\label{sec:2parcont}
\begin{figure}[t]
  \centering
  \subfigure[]{\label{sfig:2parcont_a}\includegraphics[width = 0.45 \textwidth]{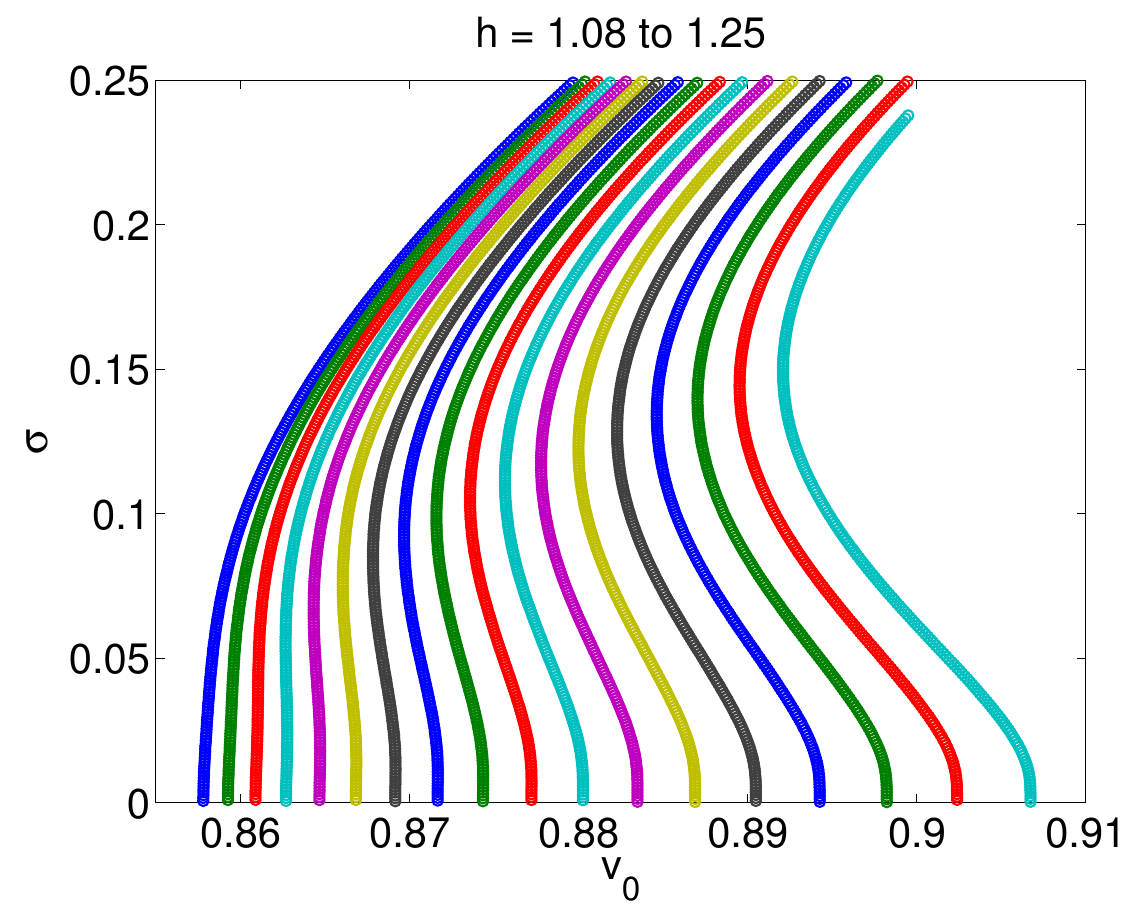}}
\hfill
  \subfigure[]{\label{sfig:2parcont_b}\includegraphics[width = 0.45 \textwidth]{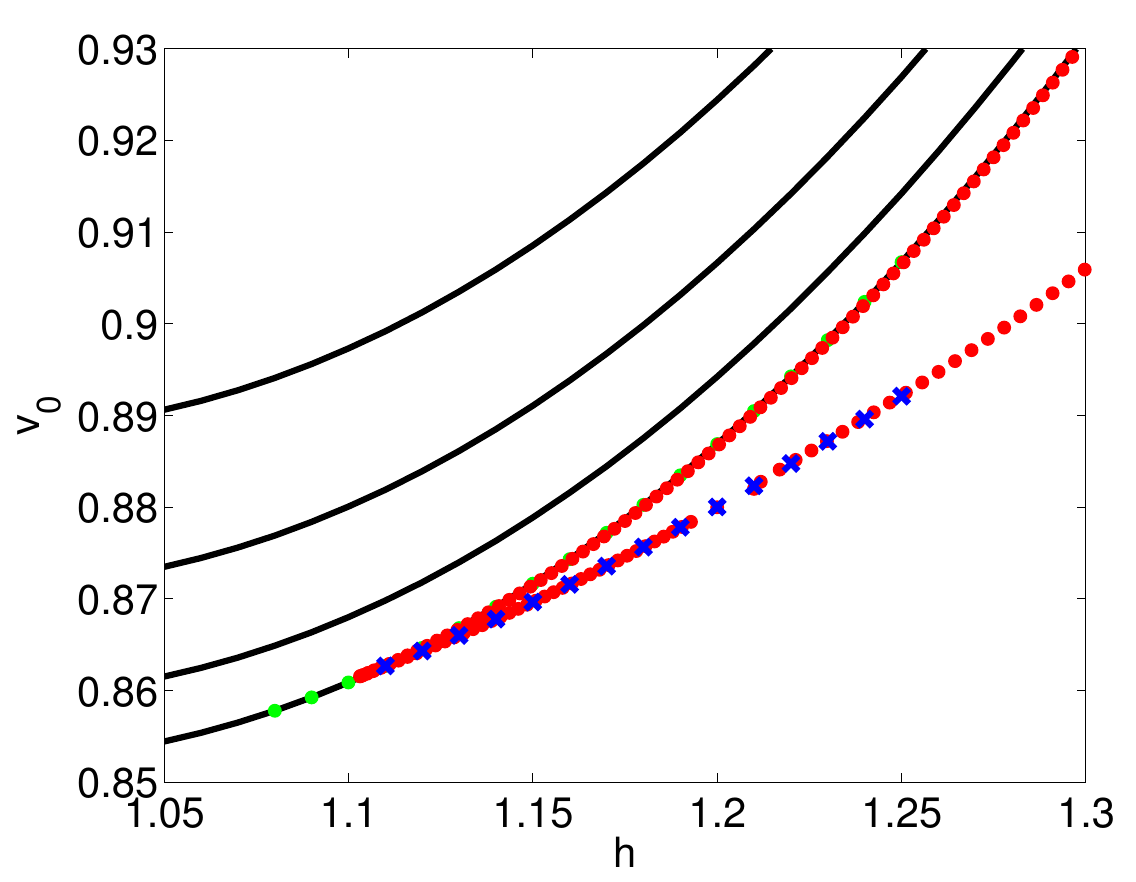}}
  \caption{\textup{(}a\textup{)} Bifurcation diagrams for $h \in
    [1.08, 1.25]$, where $h$ increases from the left curve to the
    right curve. \textup{(}b\textup{)} Two-parameter continuation of
    the fold point in the parameters $v_0$ and $h$. The blue crosses
    mark the points determined from the bifurcation diagrams of the
    one-parameter continuation, and the red circles are the results
    from a two-parameter continuation. The black lines show analytical
    results for the Hopf bifurcation at $\sigma =0$
    (cf. \eqref{eq:v0_analy}). The numerical results for the Hopf
    bifurcation (zeros from bifurcation diagrams in the left panel)
    are denoted as green dots; they are in perfect agreement with
    the analytical results. Note that in the parameter plane
    projection, the difference between the Hopf and the fold point is
    barely visible, and the first Hopf curve and the numerical data are
    obscured by the numerical data for the fold continuation.}
  \label{fig:2parcont}
\end{figure}

A two-parameter scan, showing one-parameter bifurcation diagrams in
the velocity parameter $v_0$ for different values of the safety
distance $h$, is presented in Figure
\ref{fig:2parcont}(a). The curve of folds as a result of two-parameter
continuation in Figure \ref{fig:2parcont}(a) shows how the fold merges
with another saddle-node point in a cusp. 
The system of equations for continuation of the fold is \cite{Kuznetsov2004}
\begin{equation}
  \label{eq:foldcont}
  \begin{split}
    F(\sigma,v_0,h) &= 0 \\
    F_\sigma(\sigma,v_0,h) &= 0 \\
    w^{(\sigma)} (\sigma - \hat\sigma) + w^{(v_0)} (v_0 - \hat{v}_0)
    +w^{(h)} (h - \hat{h}) &= 0
  \end{split}
\end{equation}
with the Jacobian
\begin{equation}
  \label{eq:jacobian2}
  J =
  \begin{pmatrix}
    F_\sigma & F_{v_0} & F_h\\
    F_{\sigma\sigma} & F_{v_0\sigma} & F_{h\sigma}\\
    w^{(\sigma)} & w^{(v_0)} & w^{(h)}
  \end{pmatrix}.
\end{equation}
Since derivatives of second order are needed, we apply an
approximation of second-order accuracy for the derivatives, \ie
centered differences for the parameter derivatives in $v_0$ and $h$
and one-sided second-order schemes for derivatives in $\sigma$.  We
use the one-sided second-order approximation for $F_\sigma$, because
$\sigma$ is nonnegative by definition. Details for the numerical
evaluation of the derivatives can be found in Appendix~\ref{sec:fd}.
\par
During the two-parameter continuation the Newton iteration used full
Newton steps ($\nu=1$ in \eqref{eq:newton}). Panel (b) of
Figure~\ref{fig:2parcont} shows the results; they are in perfect
agreement with the data obtained by a one-parameter continuation.  For
comparison we have included the Hopf bifurcation point of the full
microscopic system at $\sigma
=0$. The Hopf bifurcation is a pitchfork bifurcation at the
macroscopic level. However, since the standard deviation as
macroscopic measure is nonnegative by definition, it shows only the
nonnegative branches. The
analytic expression for the Hopf bifurcation parameter can be found by
linearizing system \eqref{eq:1odesys} around the uniform flow and
using the ansatz $(x_n(t), y_n(t)) = (x_n(0) \exp(i\omega t), y_n(0)
\exp(i\omega t))$.  This results in the system
\begin{eqnarray}
  \label{eq:hopf}
  i \omega x_n &=& y_n\\
  i \omega y_n &=& \tau^{-1} \left[V'\left(\frac{L}{N}\right)(x_{n+1}-x_n) - y_n\right], 
\end{eqnarray}
where $\omega$ is the frequency and $V'(\frac{L}{N})$ the first derivative of
the optimal velocity function at equilibrium.  Eliminating $x_n$ and
using the periodic boundary conditions results in
\begin{equation}
  \label{eq:hopfcond}
  \left( 1 - \frac{\omega^2 \tau}{V'\left(\frac{L}{N}\right)} + \frac{i \omega}{V'\left(\frac{L}{N}\right)} \right)^N = 1.
\end{equation}
This implicitly defines $v_0$ as a function of $h$ (through $V$) and
can be solved for our specific choice of $V$ (see \eqref{eq:ovfunc})
to yield
\begin{eqnarray}
  \label{eq:v0_analy}
  v_0 = \frac{1-\cos(2\pi j/N)}{\tau \sin^2\left(2 \pi j/N\right) \left(1- \tanh^2\left(h-\frac{L}{N}\right)\right)},
\end{eqnarray}
where $j = \lbrace 1,2,\ldots N-1 \rbrace$.  The Hopf curves for the
first four spatial frequencies ($j=1,2,3,4$) are shown in Figure
\ref{sfig:2parcont_b}. The analytical results for the first Hopf curve
are in perfect agreement with the numerical data.  Note that the
curves for the Hopf bifurcation point and the fold point are close to
each other in the parameter plane shown in
Figure~\ref{fig:2parcont}(b).

\section{Conclusion and Outlook}
\label{sec:discussion}
In this paper we have derived an implicit method for equation-free analysis
and proved its convergence for slow-fast systems with transversally
stable slow manifolds. We gave a demonstration by performing an
equation-free bifurcation analysis on a one-dimensional macroscopic
description emerging from a microscopic traffic model based on a
deterministic optimal velocity model for individual drivers. We
demonstrated that the obtained bifurcation diagrams are independent of
the lifting operator and the healing time in a suitable region. 
The bifurcation diagram shows a saddle-node bifurcation, which is
continued in a two-parameter equation-free pseudoarclength
continuation.
 Since the Hopf bifurcation, \ie the macroscopic pitch fork, is known
analytically, 
this traffic model is an ideal test case for comparison with new
numerical methods. 
The stability in Figure~\ref{fig:profiles} changes at $(v_0,\sigma) =
(0.887,0)$, \ie sign change of the eigenvalue, indicating a bifurcation.
In general, a sufficient
characterization would require checking higher-order derivatives of
the macroscopic right-hand-side $F$, which can be numerically
demanding in an equation-free computation. A detailed study of the
application of the presented implicit equation-free methods to study
pitch-fork bifurcations is a possible research direction for future
work.
\par
The proof of convergence for the implicit coarse-level time stepper
assumes that the slow manifold is transversally stable. The review
\cite{Givon2004} lists the senses in which a fast high-dimensional chaotic
or stochastic system converging in the mean can be viewed as a
slow-fast system converging to its slow manifold. In practical
applications the result from Section \ref{thm:conv} may be used as a
plausibility check: the equation-free methodology of Kevrekidis
\emph{et al} appeals to the notions of singular perturbation theory
(cf.\ the illustrative example in \cite{Kevrekidis2010}). For any
particular system under study, one can check whether this intuition is
indeed justified by testing whether the results for the implicit time
stepper given by \eqref{eq:phidef} are indeed independent of the
lifting $\lift$ and the healing time $\tskip$ if one varies both
gradually. For example, Barkley, Kevrekidis and Stuart \cite{Barkley2006} show that moment maps for simple stochastic or
  chaotic systems violate this principle in certain regions of their
  phase space.

\par
For the traffic problems studied in our paper, one long-standing
problem is the motion of several phantom jams, \ie multipulse
solutions, relative to each other. For a large number of cars
(including the $N=60$ cars we used) this motion is very slow and
therefore near impossible to observe in direct numerical simulations
(a phenomenon that is called meta-stability). An open question is
whether one can derive a computable criterion that predicts, for a
given configuration of several jams and given driver parameters, which
of those will collapse or merge and when. This criterion might be
based on the shape of the traveling wave. One particularly appealing
feature of equation-free analysis is that one can continue macroscopic
equilibria in $N$, the number of cars, using the microscopic model.
The complexity of the implicit scheme is independent of $N$. The
  increase of computational time is determined by the cost of the
  microscopic simulation with increasing $N$, since each function
  evaluation will be more costly (in our case, proportional to
  $N$). Hence, the computational complexity of the overall scheme is
  proportional to $N$.
\par
Models closer to situations of practical interest, say with more
realistic optimal velocity functions, randomly assigned driver
behavior parameters, an element of randomness in the driver behavior,
or multiple lanes, as discussed in the literature
\cite{Helbing2001,Nagatani2002,Orosz2010}, are also amenable to
equation-free analysis. This should provide additional information to
help match parameters of macroscopic models to microscopic driver
and road parameters.

\section*{Acknowledgments}
J. Starke and R. Berkemer thank Toyota CRDL for financial
support. J. Starke would also like to thank the Danish Research
Council FTP under the project number 09-065890/FTP and the Villum
Fonden under the VKR-Centre of Excellence 'Ocean Life' for financial support.
The research of J.~Sieber is supported by EPSRC grant EP/J010820/1.

\appendix

\section{Proof of Theorem \ref{thm:conv}}
\label{sec:app:conv}
For the proof of Theorem~\ref{thm:conv} we have to analyze the two
equations (for $y$ and $y_*$, respectively)
\begin{align}\label{eq:app:perturbed}
    \restrict(M_\epsilon(\tskip;\lift(y)))&=\restrict(M_\epsilon(\tskip+\deltat;\lift(x)))\mbox{,}\\
    \label{eq:app:exact}
    \restrict(M_\epsilon(\tskip;g_\epsilon(\lift(y_*))))&=
    \restrict(M_\epsilon(\tskip+\deltat;g_\epsilon(\lift(x))))\mbox{.}
\end{align}
In both equations $x\in\R^d$ enters as a parameter. For
\eqref{eq:app:exact} we have established already in Section
\ref{sec:thm} that there exists a solution $y_*$, and that it is
locally unique. Equations \eqref{eq:phistardef} and
\eqref{eq:splittimestep} gave a procedure for picking $y_*$ in a globally
unique way by starting with $y_*=x$ for $\deltat=0$ and then extending
the solution for varying $\deltat$ until one reaches the desired value
of $\deltat$. This procedure achieves unique solvability for $y_*$ for
all $x\in\dom\lift$ and for $\tskip\in[t_0,\tup/\epsilon)$ and
$\deltat\geq 0$ satisfying $\tskip+\deltat<\tup/\epsilon$. For
equation~\eqref{eq:app:perturbed} we have to prove the existence of a
solution $y$, and prove that it is close to $y_*$ (including all
derivatives with respect to $x$ up to order $k$).

In order to do this, we need to make the consequences of Fenichel's
Theorem more explicit. The Fenichel result \eqref{eq:mcebound} implies
that the map (${\cal U}$ is the neighborhood of ${\cal C}_0$, which
also contains ${\cal C}_\epsilon$)
\begin{displaymath}
  {\cal F}_\epsilon: \R\times{\cal U}\ni (\tau,u)\mapsto 
  M_\epsilon(\tau/\epsilon;g_\epsilon(u))
\end{displaymath}
is well defined and $k$ times differentiable for all $u\in{\cal U}$
and all $\tau\in\R$. This map ${\cal F}_\epsilon$ is the flow map when
restricted to the slow manifold ${\cal C}_\epsilon$, and projects all
points in the neighborhood of the slow manifold ${\cal C}_\epsilon$
along the stable fibers using $g_\epsilon$. Note that $\tau$ is the
time on the \emph{slow} time scale as we divide by $\epsilon$ in the
evaluation of the map. The derivatives of ${\cal F}_\epsilon$ with
respect to its second argument $u$ are uniformly bounded for all
$u\in{\cal U}$ as long as as $\tau\in[0,\tup]$:
\begin{equation}
  \label{eq:app:calF:bound}
  \|\partial_2^j{\cal F}_\epsilon(\tau;\cdot)\|\leq C
  \mbox{,\quad($j=0\ldots,k$, and  $\tau\in[0,\tup]$).}
\end{equation}
Correspondingly, the map
\begin{equation}
  \label{eq:app:adef}
  A_\epsilon(\tau;\cdot)=\restrict({\cal F}_\epsilon(\tau;\lift(\cdot)): 
  \dom\lift\ni x\mapsto 
  \restrict(M_\epsilon(\tau/\epsilon;g_\epsilon(\lift(x))))
\end{equation}
is well defined for all $\tau\in\R$ and locally invertible for all
$\epsilon\in[0,\epsilon_0)$ and $\tau$ satisfying $|\tau|< \tup$. Note
that the range of admissible $\epsilon$ includes $\epsilon=0$, because
the limit of the right-hand side of \eqref{eq:app:adef} for
$\epsilon=0$ is well defined as the solution of a
differential-algebraic equation on ${\cal C}_0$ on the slow time
scale. The norms of the derivatives of $A_\epsilon$ and its (locally
unique) inverse can be bounded by a uniform constant $C$ independent
of $\epsilon\in[0,\epsilon_0)$ and $\tau$ as long as $|\tau|\leq
\tup$:
\begin{align}\label{eq:app:tbound}
  \|\partial_2^jA_\epsilon(\tau;\cdot)\|&\leq C\mbox{,}&
  \|\partial_2^jA^{-1}_\epsilon(\tau;\cdot)\|&\leq C\mbox{.}
\end{align}
Similarly, the motion transversal to the slow manifold ${\cal
  C}_\epsilon$ consists of a fast decay and a slow tracking of the
dynamics on ${\cal C}_\epsilon$. Let $K<K_0$ be a given contraction
rate, and choose the upper bound $\epsilon_0$ such that the
contraction property \eqref{eq:gepsconv} of the stable fiber
projection $g_\epsilon$ holds for all $\epsilon<\epsilon_0$ and all
$u\in{\cal U}$. Then we can express the transversal component of the
flow starting from an arbitrary $u\in{\cal U}$ and $t\geq0$ in the form
\begin{equation}\label{eq:app:fiber}
  M_\epsilon(t;u)=M_\epsilon(t;g_\epsilon(u)) +
  \exp(-Kt)M^\bot_\epsilon(t;u)
\end{equation}
(this defines $M_\epsilon^\bot$).  In the right-hand side of
\eqref{eq:app:fiber} the map $M^\bot_\epsilon$ is $k$ times
differentiable with respect to its argument $u$ for all
$\epsilon\in[0,\epsilon_0)$ (including $\epsilon=0$), and the norms of
$M^\bot_\epsilon(t;u)$ and its partial derivatives
$\partial^jM^\bot_\epsilon(t;u)$ are uniformly bounded for all
$t\in[0,\infty)$, $\epsilon\in[0,\epsilon_0)$ and $u\in{\cal U}$:
\begin{align}
  \label{eq:app:mbotbound}
    \|M^\bot_\epsilon(t;x)\|&\leq C\mbox{,}&
  \|\partial_2^jM^\bot_\epsilon(t;x)\|&\leq C
  \mbox{.}
\end{align}
The prefactor $\exp(-Kt)$ can also be extracted if the smooth
restriction map $\restrict$ is applied to both terms on the left-hand
side of \eqref{eq:app:fiber}, and if we insert $\lift(x)$ for $u$. Thus,
\begin{align}
    &\restrict(M_\epsilon(t;\lift(x)))-
      \restrict(M_\epsilon(t;g_\epsilon(\lift(x))))\label{eq:app:rderiv:lhs}\\
    &
    \begin{aligned}
      =\int_0^1\partial\restrict\Big(&M_\epsilon(t;g_\epsilon(\lift(x)))
      +\rho[M_\epsilon(t;\lift(x))-M_\epsilon(t;g_\epsilon(\lift(x)))]\Big)\d \rho\\
      & \times\left[M_\epsilon(t;\lift(x))-
        M_\epsilon(t;g_\epsilon(\lift(x)))\right]
    \end{aligned}
    \label{eq:app:rderiv0}\\
    &
      =\int_0^1\partial\restrict\Big({\cal F}_\epsilon(\epsilon t;\lift(x))+
      \rho\exp(-Kt)M^\bot_\epsilon(t;\lift(x))\Big)
      \d \rho
       \exp(-Kt)
        M^\bot_\epsilon(t;\lift(x))\mbox{.}\label{eq:app:rderiv}
\end{align}
We applied the mean-value theorem to equate \eqref{eq:app:rderiv:lhs}
and \eqref{eq:app:rderiv0}. To get to the right-hand side of
\eqref{eq:app:rderiv}, we inserted the representation
\eqref{eq:app:fiber} and used the definition of the map ${\cal
  F}_\epsilon$. This right-hand side in \eqref{eq:app:rderiv} has the form
\begin{equation}\label{eq:app:rdef}
  \mbox{right-hand side of \eqref{eq:app:rderiv}}
  =\exp(-Kt)r_\epsilon(\epsilon t,t;x)\mbox{,}
\end{equation}
where the first argument of $r_\epsilon$ refers to the time dependence
of $A_\epsilon$ in the argument of $\partial\restrict$.  Note that we
have introduced the slow time scale as an additional argument into
$r_\epsilon$. We will consider $r_\epsilon(\tau,t;x)$ for arbitrary
$\tau\in[0,\tup]$ and $t\in[0,\infty)$ below, and later insert
$\tau=\epsilon t$ as a particular case.  The map
$r_\epsilon(\tau,t;x)$ is $k$ times continuously differentiable with
respect to $x$. The norm of $r_\epsilon$ and the norm of its
derivatives with respect to $x$ are uniformly bounded for
$x\in\dom\lift$, $\epsilon\in[0,\epsilon_0)$, $\tau\in[0,\tup]$, and
$t\in[0,\infty)$ because all of its ingredients have bounded
derivatives (listed in \eqref{eq:app:calF:bound},
\eqref{eq:app:mbotbound}):
\begin{align}
  \label{eq:app:rbound}
    \|r_\epsilon(\tau,t;x)\|&\leq C\mbox{,}&
  \|\partial_3^jr_\epsilon(\tau,t;x)\|&\leq C
\end{align}
for $j\in\{1,\ldots,k\}$.  Let us define the times corresponding to
$\tskip$ and $\deltat$ on the slow time scale as:
\begin{align}\label{eq:app:timescale}
  \tauskip&=\epsilon\tskip\mbox{,} &
  \Deltat&=\epsilon\deltat\mbox{.}
\end{align}
If $\tauskip$ and $\tauskip+\Deltat$ are in $[0,\tup]$, then the solution
$y_*$ of the exact flow satisfies (using the locally invertible map
$A_\epsilon$ defined in \eqref{eq:app:adef})
\begin{equation}\label{eq:app:exact:a}
A_\epsilon(\tauskip;y_*)=A_\epsilon(\tauskip+\Deltat;x)\mbox{.}
\end{equation}
Using $r_\epsilon$ and $A_\epsilon$, equation
\eqref{eq:app:perturbed} can be rewritten as
\begin{equation}\label{eq:app:perturbed:triangle}
  A_\epsilon(\tauskip;y)+s_1r_\epsilon(\tauskip,t_1;y)=
  A_\epsilon(\tauskip+\Deltat;x)+s_2r_\epsilon(\tauskip+\Deltat,t_2;x)\mbox{,}
\end{equation}
where
\begin{align}\label{eq:app:particular}
  \begin{aligned}
    s_1&=\exp(-K\tskip)\mbox{,} &
    s_2&=\exp(-K(\tskip+\deltat))\mbox{,}\\
    t_1&=\tskip\geq0\mbox{,} &
    t_2&=\tskip+\deltat\geq0\mbox{.}
  \end{aligned}
\end{align}
We will first consider solvability of
\eqref{eq:app:perturbed:triangle} with respect to $y$ for general
$s_1$ and $s_2$ close to $0$, and $t_1$, $t_2\in[0,\infty)$. This
solution $y$ will depend on the parameters $s_1$, $s_2$, $t_1$, and
$t_2$ (among others). Whenever we subsequently insert the particular
values from \eqref{eq:app:timescale} and \eqref{eq:app:particular} for
$\tauskip$, $\Deltat$, $s_1$, $s_2$, $t_1$, and $t_2$, the solution $y$
of \eqref{eq:app:perturbed:triangle} becomes also a solution of
\eqref{eq:app:perturbed}. For each of the terms, $A_\epsilon$,
$A_\epsilon^{-1}$, and $r_\epsilon$, we have uniform upper bounds
(\eqref{eq:app:tbound} and\eqref{eq:app:rbound}) for their norms and
all derivatives up to order $k$ for the entire range of arguments:
$x$, $y\in\dom\lift$, $\tauskip\in[0,\tup]$,
$\tauskip+\Deltat\in[0,\tup]$, $t_1$, $t_2\in[0,\infty)$,
and $\epsilon\in[0,\epsilon_0)$ (where $\epsilon_0$ is determined by the
choice of decay rate $K$ as given by Fenichel's Theorem). Thus, we can
use \eqref{eq:app:exact:a} and \eqref{eq:app:perturbed:triangle} to
establish the existence of $y$ and its distance to $y_*$ using the
implicit function theorem at the point $s_1=s_2=0$.

The exact solution $y_*$ is a uniformly regular solution of
\eqref{eq:app:perturbed:triangle} for $s_1=s_2=0$, all
$x\in\dom\lift$, $\epsilon\in[0,\epsilon_0)$, $\tauskip\in[0,\tup]$,
and $\Deltat\in[-\tauskip,\tup-\tauskip]$.  Thus, for small $s_1$ and
$s_2$, equation \eqref{eq:app:perturbed:triangle} has a locally unique solution
$y\in\dom\lift$ which depends smoothly on all parameters (we write
$y(x,s_1,s_2)$ to emphasize the dependence on $(s_1,s_2)\in\R^2$) such that
\begin{displaymath}
  \|\partial_1^jy(x,s_1,s_2)-\partial^j y_*(x)\|_\infty\leq C\|(s_1,s_2)\|_\infty
\end{displaymath}
($j\in\{1,\ldots,k\}$) for some constant $C$ and all $s_1$,
$s_2\in(-\rho,\rho)$ for some $\rho>0$. Consequently, if we choose
$t_0$ such that $\exp(-Kt_0)<\rho$ and decrease $\epsilon_0$ such that
$t_0<\tup/\epsilon_0$, then we have for all
$\epsilon\in(0,\epsilon_0)$, $\tskip\in[t_0,\tup/\epsilon]$,
$\deltat\in[0,\tup/\epsilon-\tskip]$, and $x\in\dom\lift$ that
\begin{align*}
  \lefteqn{\left\|\partial_1^jy\bigl(x,\exp(-K\tskip),
      \exp(-K(\tskip+\deltat))\bigr)-
      \partial^j y_*(x)\right\|_\infty}\\
  &&&\leq C\|(\exp(-K\tskip),\exp(-K(\tskip+\deltat)))\|_\infty\\
  &&&  \leq C\exp(-K\tskip)
\end{align*}
for all $j\in\{1,\ldots,k\}$. This establishes the convergence claim
of Theorem~\ref{thm:conv} since $y$ is the solution of
\eqref{eq:app:perturbed} if
$s_1=\exp(-K\tauskip/\epsilon)=\exp(-K\tskip)$,
$s_2=\exp(-K(\tauskip+\Deltat)/\epsilon)=\exp(-K(\tskip+\deltat))$,
$t_1=\tskip$, $t_2=\tskip+\deltat$, $\tauskip=\epsilon\tskip$, and
$\Deltat=\epsilon\deltat$.

\section{Parameters}
\label{sec:params}
The parameters used for the simulations are listed in Table \ref{tab:params}.
\begin{table}[ht]
  \centering
  \begin{tabular}{|c|c|}\hline
Parameter & Value/range \\ \hline
$\tau^{-1}$ & 1.7\\
$L$ & 60\\
$N$ & 60\\
$\mu$ & 0.1 \\ 
$s$ & 0.001 \\
$\delta$ & 2000\\
$\Delta t$ & -5000\\
$\tskip$ & 300\\
$v_0^*$ & 0.8, $\ldots$ ,1.0 \\
$h^*$ & 1.0, $\ldots$ ,1.7\\ \hline
  \end{tabular}
  \caption{\label{tab:params} \noindent Parameters for numerical studies. The quantities marked
    with an asterisk \textup{(}$^*$\textup{)} are bifurcation
    parameters, where the range used is noted}
\end{table}

\section{Finite Differences}
\label{sec:fd}
For the scheme \eqref{eq:foldcont},
$F$ is evaluated at the 17 points
\begin{equation}
\begin{aligned}
  \label{eq:evpoints}
&  1:(\sigma, v_0,h),& \qquad \\  
  &2:(\sigma+\Delta \sigma, v_0,h), \qquad
  &3&:(\sigma+2\Delta \sigma, v_0,h), \\
&  4:(\sigma+3\Delta \sigma, v_0,h), \qquad
  &5&:(\sigma+4\Delta \sigma, v_0,h), \\
&  6:(\sigma, v_0-\Delta v_0,h), \qquad   
  &7&:(\sigma, v_0+\Delta v_0,h), \\
&  8:(\sigma+\Delta \sigma, v_0-\Delta v_0,h), \qquad
  &9&:(\sigma+\Delta \sigma, v_0+\Delta v_0,h), \\
&  10:(\sigma+2\Delta \sigma, v_0-\Delta v_0,h), \qquad
  &11&:(\sigma+2\Delta \sigma, v_0+\Delta v_0,h),\\
&  12:(\sigma, v_0,h-\Delta h), \qquad   
  &13&:(\sigma, v_0,h+\Delta h), \\
&  14:(\sigma+\Delta \sigma, v_0,h-\Delta h), \qquad
  &15&:(\sigma+\Delta \sigma, v_0,h+\Delta h), \\
&  16:(\sigma+2\Delta \sigma, v_0,h-\Delta h), \qquad
  &17&:(\sigma+2\Delta \sigma, v_0,h+\Delta h).
  \end{aligned}
\end{equation}
where $\Delta \sigma = \Delta v_0 = \Delta h
= 0.001$ are offsets for the approximation.  One can use the following
second-order accuracy scheme to compute the derivatives (for better
readability, the points are just referred to by their number, \eg
$F_7=F(\sigma, v_0+\Delta v_0,h)$):
\begin{eqnarray*}
  \label{eq:compjac}
  F_\sigma &=& \frac{-3F_1 +4 F_2- F_3}{2\Delta \sigma} \\
  F_{v_0} &=& \frac{F_7-F_6}{2 \Delta v_0}\\
  F_{h} &=& \frac{F_{13}-F_{12}}{2 \Delta h}\\
  F_{\sigma\sigma} &=& \frac{-3(-3F_1 + 4F_2 -F_3) +4(-3F_2 +
    4F_3 -F_4)-(-3F_3+4F_4-F_5)}{4(\Delta \sigma)^2} \\
  F_{v_0\sigma} &=& \frac{(-3F_7 + 4F_9
    -F_11)-(-3F_6+4F_8-F_{10})}{4\Delta \sigma \Delta v_0} \\
  F_{h\sigma} &=& \frac{(-3F_{13} + 4F_{15} -F_{17})-(-3F_{12}+4F_{14}-F_{16})}{4\Delta \sigma \Delta h}
\end{eqnarray*}


\bibliographystyle{siam}
\bibliography{traffic_implicit}

\end{document}